\documentclass[a4paper,11pt,leqno,twoside]{amsart}
\usepackage{honoriro}
\begin{document}%
\section{Introduction}%
\subsection{Statement of the results}%
The correspondence principle in quantum mechanics suggests a way to study a classical system via its semi-classical limit of quantization. For instance, let $X$ be a compact Riemannian manifold. We can choose an orthonormal basis $\left(f_j\right)_{j\geq 0}$ of $L^2(X)$ satisfying
\begin{equation*}
\forall j\geq 0, \quad \Delta(f_j)=\lambda_j f_j.
\end{equation*}
where $\Delta$ is the Laplace-Beltrami operator on $X$ and $0=\lambda_0<\lambda_1\leq\lambda_2\leq\dots$ is its spectrum. If $G^t$ is the geodesic flow on $X$ then its quantization is $-h^2\Delta$, where $h$ is Planck's constant. Thus it is very natural to attempt to understand the asymptotic behaviour of the eigenfunctions of $\Delta$. 
\par
A classical question here -- suggested by the correspondence principle -- is to bound $\abs{\abs{f_j}}_\infty$ as $\lambda_j\to \infty$. (See \cite{MR1886720} and \cite{MR1321639} for more details.) A.~Seeger and C.~Sogge proved in \cite{MR1017329} a very general and qualitative bound, essentially sharp, in the case of compact Riemannian surfaces.
\par
If $X$ is a compact locally symmetric space then P.~Sarnak proved in \cite{SaMo} the generic bound
\begin{equation*}
\abs{\abs{f_j}}_\infty\ll\lambda_j^{\left(\text{dim}(X)-\text{rank}(X)\right)/4}
\end{equation*}
provided $f_j$ is the joint eigenfunction of all the algebra of the invariant differential operators.
\par
In \cite{MR1324136}, H.~Iwaniec and P.~Sarnak proved a bound sharper than that of A.~Seeger and C.~Sogge for certain Hecke eigenfunctions on arithmetic surfaces which are the quotient of the upper-half plane by a congruence subgroup of $SL_2(\Z)$; they took advantage of the fact that some additional symmetries, the Hecke correspondences, act on these surfaces and one can take an orthonormal basis of Hecke eigenfunctions. The Laplace-Beltrami operator in this context is the hyperbolic
Laplacian.
\par
Following this foundational result, the sup-norm problem in the eigenvalue aspect has since been considered in various settings. For instance, S.~Koyama investigated the case of quotients of the three-dimensional hyperbolic space by arithmetic subgroups in \cite{MR1324641} and proved similar results, which have been improved by V.~Blomer, G.~Harcos and D.~Milicevic in \cite{BlHaMil}. J.~Vanderkam \cite{MR1440572} and later on V.~Blomer and P.~Michel (\cite{MR2852302})considered the case of the sphere and of the ellipsoids. S.~Marshall considered the sup-norm problem restricted to totally geodesic submanifolds in \cite{Ma1} and in \cite{Ma2}. V.~Blomer and A.~Pohl considered for the first time a manifold of higher rank and solved the case of Hecke Siegel Maass cusp form of genus $2$ for $Sp_4(\Z)$ in \cite{BlPo}. 
\par
We will focus on another non-compact Riemannian symmetric space of dimension $5$ and rank $2$, which is 
\begin{equation*}
X=SL_3(\Z)\backslash SL_3(\R)\slash SO_3(\R).
\end{equation*}
\par
In this manuscript, we provide a proof of a non-trivial explicit quantitative upper bound for a $SL_3(\Z)$ Hecke-Maass cusp form at a generic point $z$ in a fixed compact subset of $X$. These forms are Maass forms since they are eigenfunctions of the algebra of invariant differential operators and Hecke forms since we assume they are eigenfunctions under the Hecke operators. Specifically, we establish the following result.
\begin{theoint}\label{theo_mainresult}
  Let $\Phi$ be an $L^2$-normalized and tempered $SL_3(\Z)$ Hecke-Maass cusp 
  form on $X$ with Laplace eigenvalue $\lambda$ 
  and type $(\nu_1, \nu_2)$ in $i\R^2$ satisfying $|\nu_1 - \nu_2| \ll 1$.
  Let $C$ be a fixed compact in $X$. One has
\begin{equation*}
\abs{\abs{\Phi_{\vert C}}}_\infty\ll_{C,\epsilon}\lambda^{(5-2)/4-1/76+\epsilon}
\end{equation*}
for all $\epsilon>0$. 
\end{theoint}
\par
Several works related to this problem have to be mentioned. In \cite{MR3384442}, V.~Blomer and P.~Maga proved a qualitative non-trivial bound for the sup-norm of $PGL(4)$ Hecke-Mass cusp forms restricted to compact sets and in \cite{MR3518551}, they proved a qualitative non-trivial bound for the sup-norm of $PGL(n)$ Hecke-Mass cusp forms restricted to compact sets for $n\geq 5$. In both works, the subconvexity exponent is not computed. Very recently, V.~Blomer, G.~Harcos and P.~Maga proved in \cite{BlHaMa} a quantitative bound for the global (namely without any restriction to compact sets) sup-norm of $GL(3)$ Hecke-Mass cusp forms. More explicitly, they proved that
\begin{equation*}
\abs{\abs{\Phi}}_\infty\ll_{\epsilon}\lambda^{(5-2)/4+9/40+\epsilon}
\end{equation*}
for any $L^2$-normalized and tempered $SL_3(\Z)$ Hecke-Maass cusp form $\Phi$ on $X$ with Laplace eigenvalue $\lambda$ and for any $\epsilon>0$.
\par
The method of proof builds on generalizations of the work of H.~Iwaniec and P.~Sarnak in \cite{MR1324136}, i.e. one studies a smooth amplified second moment, which comes from the spectral expansion of an automorphic kernel, which itself has a geometric expansion. This is usually referred to as the pre-trace formula.
\par
An amount of time is devoted to the construction of a relevant function on the spectral side of the pre-trace formula. In particular, one has to bound its inverse Helgason transform in the different domains of the positive Weyl chamber. This relies on the spherical inversion formula and on a systematic study of the $GL(3)$ spherical function itself done by S.~Marshall in \cite{Ma3}.
\par
Finally, the geometric side of the amplified pre-trace formula is bounded thanks to a counting lemma, which is the analogue of the one seen in \cite{BlPo}.
\subsection{Organization of the paper}%
Section \ref{sec_background} contains the knowledge on Lie groups and Lie algebras required for this work and all the relevant notations. Section \ref{sec_aptf} briefly explains the strategy of the proof and states an amplified pre-trace formula. The background on the $GL(3)$ Hecke algebra is given in Section \ref{sec_Hecke}. Moreover, several linearizations of compositions of some Hecke operators, which are required to make the amplification effective and done in \cite{MR3525536}, are recalled. In Section \ref{sec_test}, the function which occurs on the spectral side of the amplified pre-trace formula is constructed and several estimates for its inverse Helgason transform are proven. Section \ref{sec_first_estimate} contains a first bound for the geometric side of the amplified pre-trace formula, based on the results done in the previous sections. The counting lemma required to complete this bound is given in Section \ref{sec_counting}. The end of the proof of Theorem \ref{theo_mainresult} appears in the final section.
\begin{notations}
The main parameters in this work are a positive real number $T$, which goes to infinity and a positive integer $L$ (a power of $T$ determined at the very final step) which goes to infinity with $T$. Thus, if $f$ and $g$ are some $\C$-valued functions on $\R^2$ then the symbols $f(T,L)\ll_{A}g(T,L)$ or equivalently $f(T,L)=O_A(g(T,L))$ mean that $\abs{f(T,L)}$ is smaller than a constant, which only depends on $A$, times $g(T,L)$. Similarly, $f(T,L)=o(1)$ means that $f(T,L)\to 0$ as $T$ goes to infinity among the positive real numbers.
\par
We will denote by $\epsilon$ a positive constant whose value may vary from one line to the next one.
\end{notations}
\begin{merci}
\par
The authors would like to thank the referee for her or his careful reading of the manuscript.
\par
They also would like to thank V.~Blomer, F.~Brumley, J.~Cogdell, \'{E}.~Fouvry, H.~Iwaniec, E.~Kowalski, E.~Lapid, S.~Marshall, P.~Michel, A.~Pohl, P.~Sarnak and R.~J.~Stanton for stimulating exchange related to this project.
\par
This paper was worked out at several places: while the first, third and fourth authors attended the workshop "Analytic theory of $GL(3)$ automorphic forms and applications" by the American Institute of Mathematics in Palo Alto, while the first, third and fourth authors were invited by Forschungsinstitut für Mathematik (FIM, ETH) in Zürich, while the first and third authors were invited by Université Blaise Pascal (Laboratoire de Mathématiques) in Clermont-Ferrand, while the third and fourth authors were invited by The Ohio State University (Department of Mathematics) in Columbus, while the second author was invited by Université de Bordeaux. We would like to thank all these institutions for their hospitality and inspiring working conditions.
\par
R.~Holowinsky was supported by the Sloan fellowship BR2011-083 and the NSF grant DMS-1068043.
\par
K.~Nowland is supported as a Graduate Research Associate by The Ohio State University (Department of Mathematics).
\par
The research of G.~Ricotta was supported by a Marie Curie Intra European Fellowship within the 7th European Community Framework Programme. The grant agreement number of this project, whose acronym is ANERAUTOHI, is PIEF-GA-2009-25271. He would like to thank ETH and its entire staff for the excellent working conditions.
\par
The third and fourth authors are financed by the ANR Project Flair ANR-17-CE40-0012. 
\par
This material is partially based upon research supported by the Chateaubriand Fellowship of the Office for Science \& Technology of the Embassy of France in the United States.
\end{merci}
\section{Background on Lie groups and Lie algebras}\label{sec_background}%
Let $G:= SL_3(\R)$ and
\begin{equation*}
A=\left\{a=\begin{pmatrix}
a_1 & & \\
& a_2 & \\
& & a_3
\end{pmatrix}\in M_3(\R), \text{det}(a)=1, \forall i\in\{1,2,3\}, a_i>0\right\},
\end{equation*}
whose Lie algebra is
\begin{equation*}
\mathfrak{a}=\left\{H=\begin{pmatrix}
h_1 & & \\
& h_2 & \\
& & h_3
\end{pmatrix}\in M_3(\R), \text{Tr}(H)=0\right\},
\end{equation*}
whose complexification is denoted by $\mathfrak{a}_\C$. Let
\begin{equation*}
N=\left\{n=\begin{pmatrix}
1 & x_1 & x_3\\
& 1 & x_2 \\
& & 1
\end{pmatrix}\in M_3(\R)\right\}
\end{equation*}
and $K:=SO_3(\R)$ be one of the maximal compact subgroups of $G$.
\par
The \emph{Iwasawa decomposition} of $G$ is given by $G=NAK$. If $g=nak$ then one denotes by
\begin{equation*}
\IwK(g)=k \text{ and } \IwA(g)=a.
\end{equation*}
\par
The set
\begin{equation*}
\beta\coloneqq\left\{H_{1,2}=E_{1,1}-E_{2,2},H_{2,3}=E_{2,2}-E_{3,3}\right\}
\end{equation*} 
is a basis of the $2$-dimensional $\R$-vector space $\mathfrak{a}$ where $E_{i, j}$ the matrix with all zero entries except for a $1$ in the $i$th row and $j$th column. The \emph{Killing form}
\begin{equation*}
B(H,H')=6\text{Tr}(HH')
\end{equation*}
is a positive definite quadratic form on $\mathfrak{a}$. The same properties hold for $\mathfrak{a}_\C$, the only difference being that the Killing form is a non-degenerate bilinear symmetric form on $\mathfrak{a}_\C$.
\par
The $\R$-linear forms
\begin{equation*}
\alpha_{i,j}(H)=h_i-h_j
\end{equation*}
for $1\leq i,j\leq 3$ belong to $\mathfrak{a}^\ast$ and
\begin{equation*}
\beta_1^\ast\coloneqq\left\{\alpha_1^{+}=\alpha_{1,2},\alpha_{2}^{+}=\alpha_{2,3}\right\}
\end{equation*}
is a basis of the $2$-dimensional $\R$-vector space $\mathfrak{a}^\ast$, whose elements are called the \emph{simple positive roots}. The last positive root is $\alpha_3^{+}=\alpha_1^{+}+\alpha_2^{+}$. The multiplicative roots on $A$ are
\begin{equation*}
\forall i\in\{1,2,3\},\quad\alpha_i(a)=e^{\alpha_i^+(\log{(a)})}.
\end{equation*}
\begin{figure}[H]
\centering 
\scalebox{0.5} 
{\includegraphics{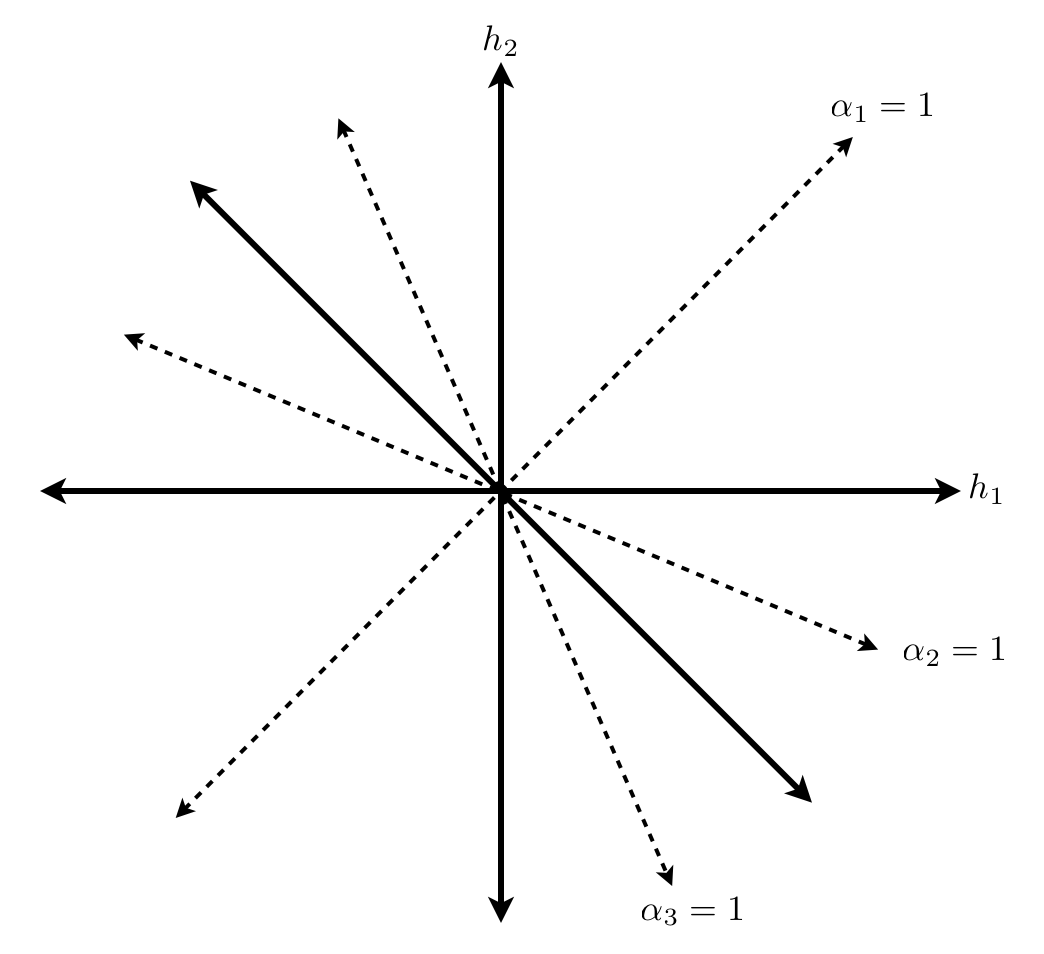}} 
\caption{$A=\exp{H}$.}
\label{fig:additive} 
\end{figure} 
Another basis of $\mathfrak{a}^\ast$ is given by
\begin{equation*}
\beta_2^\ast\coloneqq\left\{\lambda_1,\lambda_{2}\right\}
\end{equation*}
where
\begin{equation*}
\lambda_1(H)=h_1,\quad\lambda_2(H)=h_1+h_2.
\end{equation*}
One can check that $\beta_2^\ast$ is the dual basis of $\beta$. The same properties hold for $\mathfrak{a}_\C^\ast$.
\par
The Killing form being positive definite on $\mathfrak{a}$, one can identify canonically $\mathfrak{a}$ and $\mathfrak{a}^\ast$, in the sense that
\begin{equation*}
\forall\lambda\in\mathfrak{a}^\ast, \exists! H_\lambda\in\mathfrak{a},\quad\lambda=B(H_\lambda,\ast).
\end{equation*}
In addition, one can transfer the Killing form to $\mathfrak{a}^\ast$ by the formula
\begin{equation*}
\forall(\lambda,\mu)\in\left(\mathfrak{a}^\ast\right)^2,\quad B(\lambda,\mu)\coloneqq B(H_{\lambda},H_{\mu}).
\end{equation*}
The basis $\{6\lambda_1, 6\lambda_2\}$ is the $B$-dual basis of the basis $\beta_{1}^\ast$, in the sense that
\begin{equation*}
B\left(6\lambda_i,\alpha_j^+\right)=\delta_{i,j}
\end{equation*}
for $1\leq i,j\leq 2$. The same properties hold for $\mathfrak{a}_\C^\ast$ since the Killing form is non-degenerate on $\mathfrak{a}_\C^\ast$.
\par
One can also define a positive definite quadratic form on $\mathfrak{a}_\C^\ast$ as follows. Obviously,
\begin{equation*}
\mathfrak{a}_\C^\ast=\mathfrak{a}^\ast\oplus i\mathfrak{a}^\ast.
\end{equation*}
If $\lambda=\lambda_\R+\lambda_I$ with respect to this decomposition then the conjugate of $\lambda$ is defined to be $\lambda^{\text{conj}}=\lambda_\R-\lambda_I$. The bilinear symmetric form on $\mathfrak{a}_\C^\ast$ given by
\begin{equation*}
f\left(\lambda_1,\lambda_2\right)=B\left(\lambda_1,\lambda_2^{\text{conj}}\right)
\end{equation*}
is positive definite and the induced norm is
\begin{equation*}
\abs{\abs{\lambda}}=\sqrt{f(\lambda,\lambda)}=\sqrt{\abs{\abs{\lambda_\R}}^2+\abs{\abs{\lambda_I}}^2}.
\end{equation*}
\par
The basis $\beta_2^\ast$ is the one that will be used to find an explicit integral representation for the spherical function. If $\lambda$ belongs to $\mathfrak{a}_\C^\ast$ then there exists a unique pair $s=(s_1,s_2)$ of complex numbers satisfying
\begin{equation}\label{eq_lambda_s}
\lambda=s_1\lambda_1+s_2\lambda_2.
\end{equation}
One writes $\lambda=\lambda_s$.
\par
One gets for free an explicit parametrization of the multiplicative characters on $A$ via the exponential map. If $a$ belongs to $A$ then one can define
\begin{equation*}
p_1(a)=a_1,\quad p_2(a)=a_1a_2.
\end{equation*}
For $s=(s_1,s_2)$ a pair of complex numbers, the Selberg character of parameter $s$ is given by
\begin{equation*}
p_s(a)=p_1(a)^{s_1}p_2(a)^{s_2}.
\end{equation*}
A famous one is the module given by $\delta=p_{(1,1)}^2$. All the multiplicative characters of $A$ are of this shape. If $\chi:A\to\C$ is a multiplicative character then there exists a unique pair $s=(s_1,s_2)$ of complex numbers satisfying
\begin{equation*}
\chi=p_s.
\end{equation*}
Note that
\begin{equation*}
\exp\circ\lambda_s=p_s\circ\exp.
\end{equation*}
\par
The \emph{Weyl group} $W$ of $G$ is the quotient of the normalizer of $A$ in $K$ by the centralizer of $A$ in $K$. Its action on $A$ can be identified with the action of the symmetric group $\sigma_3$ by permuting the diagonal elements of the diagonal matrices in $A$. $W$ also acts on $\mathfrak{a}$ and $\mathfrak{a}_\C$ by permuting the diagonal matrices of these vector spaces. A fundamental domain for this action of $W$ on $A$ is given by the \emph{positive Weyl chamber}
\begin{equation*}
A_+\coloneqq\left\{a\in A, \alpha_1(a)>1, \alpha_2(a)>1\right\}.
\end{equation*}
This action is transferred to an action on the group of multiplicative characters on $A$ as follows. Recall that for $s\in\C^2$, the multiplicative character $\chi_s$ can be identified with the $\C$-linear function $\lambda_s$. For $w\in W$, one can define the $\C$-linear function $w.\lambda_s$ by
\begin{equation*}
w.\lambda_s=B\left(w.H_{\lambda_s},\ast\right).
\end{equation*}
In other words, $H_{w.\lambda_s}=w.H_{\lambda_s}$. The multiplicative character $w.p_s$ is the multiplicative character associated to $w.\lambda_s$, namely
\begin{equation*}
\left(w.p_s\right)(a)=\exp{\left((w.\lambda_s)(\log{(a)})\right)}.
\end{equation*}
Equivalently, $W$ acts on $\C^2$ by the explicit formulas given by
\begin{eqnarray*}
(1,2).s & = & (-s_1,s_1+s_2), \\
(1,3).s & = & (-s_2,-s_1), \\
(2,3).s & = & (s_1+s_2,-s_2), \\
(1,2,3).s & = & (s_2,-s_1-s_2), \\
(1,3,2).s & = & (-s_1-s_2,s_1).
\end{eqnarray*}
\par
Recall that the \emph{Cartan decomposition} of $G$ is $G=KAK$. If $g=k_1ak_2$ then one has a simple formula for the \emph{geodesic distance} on the Riemannian manifold $G\slash K$ between $g$ and the identity matrix $I$. Since $I$ is fixed by the action of $K$, our distance function will only depend on the entries of $a$ as 
\begin{equation}\label{eq_dist}
d(g,I)^2:=\log^{2}{(a_1)}+\log^{2}{(a_2)}+\log^{2}{(a_3)}.
\end{equation}
Up to a constant, this notion of distance comes from taking 
$\theta(X)\coloneqq-{}^tX$ as a Cartan involution on the Lie algebra 
and defining a notion of size as $B(X,-\theta(X))$ with corresponding
distance between
$X$ and $Y$ as $B(X-Y,-\theta(X - Y))$.
In terms of the multiplicative roots, this becomes
\begin{equation*}
d(g,I)^2=\frac{2}{3}\left(\log^2{(\alpha_1(a))}+ \log{(\alpha_1(a))}\log{(\alpha_2(a))}+\log^2{(\alpha_2(a))}\right).
\end{equation*}
\section{The amplified pre-trace formula}\label{sec_aptf}%
Let $\Phi_{j_0}$ be our favorite $SL_3(\Z)$ Hecke-Maass cusp form of type$\nu_{j_0}=(\nu_{{j_0},1},\nu_{{j_0},2})\in\C^2$. The background on these objects can be found in \cite{MR2254662}. One can include $\Phi_{j_0}$ in an orthonormal basis of $SL_3(\Z)$ Hecke-Maass cusp forms $\left(\Phi_j\right)_{j\geq 0}$, the type of each $\Phi_j$ being $\nu_j=(\nu_{j,1},\nu_{j,2})\in\C^2$ for $j\geq 0$.
\par
Let $k$ be a smooth and compactly supported bi-$K$-invariant function on $G$ satisfying the following properties.
\begin{itemize}
\item
For $j\geq 0$, $\mathcal{H}(k)(\nu_j)\geq 0$ where $\mathcal{H}(k)$ is the Helgason transform of $k$ (see Section \ref{sec_test}).
\item
$\mathcal{H}(k)$ is non-negative on the continuous spectrum of $X$.
\item
$\mathcal{H}(k)(\nu_{j_0})\gg 1$.
\end{itemize}
\par
Let $K(z,z^\prime)$ be the automorphic kernel given by
\begin{equation} \label{automorphickernel}
K(z,z^\prime)\coloneqq\sum_{\gamma\in GL_3(\mathbb{Z})/\{\pm I\}}k(z^{-1}\gamma z^\prime)
\end{equation}
for all $z$ and $z'$ in $G$. This function is left-$SL_3(\Z)$-invariant and right-$K$-invariant with respect to each variable $z$ and $z'$.
\par
Spectrally decomposing via a pre-trace formula, one gets that 
\begin{equation}\label{eq_pre_trace}
K(z,z^\prime)=\sum_{j\geq 0}\mathcal{H}(k)(\nu_j)\Phi_j(z^\prime)\overline{\Phi_j(z)}+\dots
\end{equation}
where $\dots$ stands for the contribution of the continuous spectrum.
\par
Let $I$ be a suitable finite subset of $\mathbb{N}^2$ and let $\alpha=(\alpha_{m,n})_{(m,n)\in I}$ be a suitable sequence of complex numbers which will be chosen later. Assume the existence of linear operators $T_{m,n}$ and $T_{m,n}^\ast$ such that
\begin{eqnarray} \label{operatorconditions}
T_{m,n}(\Phi_j) & = & a_j(m,n)\Phi_j, \\
T_{m,n}^\ast(\Phi_j) & = & \overline{a_j(m,n)}\Phi_j
\end{eqnarray}
for $(m,n)\in I$. We shall later choose $a_j(m,n)$ to be the Hecke eigenvalues of certain Hecke operators $T_{m,n}$. Defining
$$
A_j(\alpha)\coloneqq\sum_{(m,n)\in I}\alpha_{m,n}a_j(m,n),
$$
one has that
\begin{multline*}
\sum_{j\geq 0}\left\vert A_j(\alpha) \right\vert^2\mathcal{H}(k)(\nu_j)\Phi_j(z^\prime)\overline{\Phi_j(z)}+\ldots=\sum_{\substack{(m_1,n_1)\in I \\
(m_2,n_2)\in I}}\alpha_{m_1,n_1}\overline{\alpha_{m_2,n_2}} 
\sum_{j\geq 0}\mathcal{H}(k)(\nu_j)a_j(m_1,n_1)\overline{a_j(m_2,n_2)}\Phi_j(z^\prime)\overline{\Phi_j(z)}+\ldots
\end{multline*}
upon expanding the square and where $\ldots$ stands for the contribution of the continuous spectrum of $X$.
\par
Fix $z$ and consider the previous equality as an equality of functions of $z^\prime$. One has
\begin{align*}
  \sum_{j\geq 0}\left\vert A_j(\alpha) \right\vert^2 & \mathcal{H}(k)(\nu_j)\Phi_j(z^\prime)\overline{\Phi_j(z)} + \ldots \\
&=\sum_{\substack{(m_1,n_1)\in I \\
(m_2,n_2)\in I}}\alpha_{m_1,n_1}\overline{\alpha_{m_2,n_2}}\sum_{j\geq 0}\mathcal{H}(k)(\nu_j)\left[\left(T_{m_2,n_2}^\ast\circ T_{m_1,n_1}\right)(\Phi_j)\right](z^\prime)\; \overline{\Phi_j(z)} + \ldots.
\end{align*} 
\par
By \eqref{eq_pre_trace}, this gives
\begin{align*}
  \sum_{j\geq 0}\left\vert A_j(\alpha) \right\vert^2 \mathcal{H}(k)(\nu_{j}) & \Phi_j(z^\prime)\overline{\Phi_j(z)}+\ldots \\
&= \sum_{\substack{(m_1,n_1)\in I \\
(m_2,n_2)\in I}}\alpha_{m_1,n_1}\overline{\alpha_{m_2,n_2}}\left[\left(T_{m_2,n_2}^\ast\circ T_{m_1,n_1}\right)(K(z,\ast))\right](z^\prime).
\end{align*}
Here we have used the fact that the Hecke operators $T_{m,n}$ act on the Eisenstein series in the continuous spectrum in the same way in which they act on Hecke-Maass cusp forms. The left-hand side of this formula is the \emph{spectral side} whereas the right-hand side is the \emph{geometric side} of the amplified pre-trace formula.
\par
Choosing $z=z^\prime$, one makes use of positivity of the summand and estimates the size of any single $\Phi_{j_0}(z)$ by the following inequality
\begin{equation}\label{eq_amplified_pre_trace}
\left\vert A_{j_0}(\alpha) \right\vert^2\mathcal{H}(k)(\nu_{j_0})\left\vert\Phi_{j_0}(z)\right\vert^2\leq\sum_{\substack{(m_1,n_1)\in I \\
(m_2,n_2)\in I}}\alpha_{m_1,n_1}\overline{\alpha_{m_2,n_2}}\left[\left(T_{m_2,n_2}^\ast\circ T_{m_1,n_1}\right)(K(z,\ast))\right](z).
\end{equation}
Therefore, everything boils down to bounding the geometric side of the amplified pre-trace formula.
\par
We will choose the coefficients $\alpha_{m,n}$ such that $\left\vert A_{j_0}(\alpha)\right\vert$ is bounded below by a small power of the main parameter $T$. We will also choose the coefficients $a_j(m,n)$ such that it will be possible to linearize the composition $T_{m_2,n_2}^\ast\circ T_{m_1,n_1}$. See Section \ref{sec_Hecke} for an explicit description of all these parameters.
\par
We will not choose the function $k$ occurring in \eqref{eq_amplified_pre_trace} but instead the function\footnote{Actually, similarly to what did H.~Iwaniec and P.~Sarnak in \cite[Section 1]{MR1324136}, we will choose the inverse Fourier transform of $\mathcal{H}(k)$.} $\mathcal{H}(k)$ with the required properties and we will prove the needed estimates for the corresponding function $k$ in order to bound the geometric side of the amplified pre-trace formula.
\section{The Hecke algebra}\label{sec_Hecke}%
\subsection{Linearizations of Hecke operators}%
For $g$ a matrix in $GL_3(\mathbb{Q})$, the Hecke operator $T_g$ acts on a $\mathbb{C}$-valued function $f$ defined on $G$, which is left-$SL_3(\Z)$-invariant and right-$K$-invariant, by the formula
\begin{equation*}
\left(T_g(f)\right)(z)=\sum_{\delta\in GL_3(\Z)\setminus GL_3(\Z)gGL_3(\Z)}f\left(\frac{1}{\text{det}(\delta)^{1/3}}\delta z\right)
\end{equation*}
for all $z$ in $SL_3(\R)$. Note that on the one hand, the double coset $GL_3(\Z)gGL_3(\Z)$ is a finite union of left $GL_3(\Z)$ cosets since $g$ belongs to $GL_3(\mathbb{Q})$ and on the other hand, $T_g$ is well-defined since its definition does not depend on a choice of representatives of the quotient set because $f$ is left-$SL_3(\Z)$-invariant. The resulting new function $T_g(f)$ remains left-$SL_3(\Z)$-invariant and right-$K$-invariant. The fact that $g$ is allowed to have rational coefficients and not only integer ones is required for the theory since the adjoint with respect to the Petersson inner product of $T_g$ is $T_{g^{-1}}$.
\par
One can compute the action of such Hecke operator $T_g$ on the automorphic kernel as follows. Let us fix a matrix $z$ in $G$. One successively gets
\begin{align}\label{eq_auto_action}
\left(T_g(K(z,\ast))\right)(z^\prime) & =\sum_{\delta\in GL_3(\Z)\setminus GL_3(\Z)gGL_3(\Z)}\sum_{\gamma\in GL_3(\Z)/\{\pm I\}}k\left(\frac{1}{\text{det}(\delta)^{1/3}}z^{-1}\gamma\delta z^\prime\right) \\
& =\sum_{\delta\in GL_3(\Z)\setminus GL_3(\Z)gGL_3(\Z)}\sum_{\gamma\in GL_3(\Z)/\{\pm I\}}k\left(\frac{1}{\text{det}(\gamma\delta)^{1/3}}z^{-1}\gamma\delta z^\prime\right) \\
& =\sum_{\rho\in GL_3(\Z)gGL_3(\Z)/\{\pm I\}}k\left(\frac{1}{\text{det}(\rho)^{1/3}}z^{-1}\rho z^\prime\right)\label{eq_K_Hecke}
\end{align}
for each matrix $z'$ in $G$. The equation \eqref{eq_K_Hecke} reveals that we should have a clear understanding of the double coset of $g$.
\par
The main reference is \cite{MR0340283}. Let $g=\left[g_{i,j}\right]_{1\leq i,j\leq 3}$ be a matrix of size $3$ with integer coefficients and $k\leq 3$ be a positive integer. Let $I_{k}$ be the finite set of all $k$-tuples $\{i_1,\dots,i_k\}$ satisfying $1\leq i_1<\dots<i_k\leq 3$. If $\omega$ and $\tau$ are two elements of $I_{k}$ then $g(\omega,\tau)$ will denote the $k\times k$ determinantal minor of $g$ whose row indices are the elements of $\omega$ and whose column indices are the elements of $\tau$. The $k$-th \emph{determinantal divisor} of $g$ say $d_k(g)$ is defined by
\begin{equation*}
d_k(g)\coloneqq\begin{cases}
0 & \text{if $\forall(\omega,\tau)\in I_{k}^2, g(\omega,\tau)=0$,} \\
\ggcd\left\{g(\omega,\tau), (\omega,\tau)\in I_{k}^2\right\} & \text{otherwise}
\end{cases}
\end{equation*}
where the $\ggcd$ is chosen to be positive. In particular,
\begin{equation*}
d_1(g)=\ggcd\left\{\left\vert g_{i,j}\right\vert, 1\leq i,j\leq 3\right\},\quad d_3(A)=\left\vert\dt(g)\right\vert.
\end{equation*}
These quantities are useful since they completely determine a given double coset. More precisely, a matrix $h$ of size $3$ with integer coefficients belongs to $GL_3(\Z)gGL_3(\Z)$ if and only if
\begin{equation*}
\forall 1\leq k\leq 3,\quad d_k(h)=d_k(g).
\end{equation*}
The determinantal divisors satisfy the divisibility properties
\begin{equation}\label{eq_det_div_1}
\forall 1\leq k\leq 2,\quad d_k(A)^2\mid d_{k-1}(A)d_{k+1}(A)
\end{equation}
with the convention $d_0(A)=1$ and 
\begin{equation}\label{eq_det_div_2}
d_1(A)^k\mid d_k(A)
\end{equation}
for $1\leq k\leq 3$.
\par
For $n$ a positive integer, the $n$-th normalized Hecke operator is defined by
\begin{equation*}
T_n\coloneqq\frac{1}{n}\sum_{\substack{g=\text{diag}(y_1,y_2,y_3) \\
y_1\mid y_2\mid y_3 \\
y_1y_2y_3=n}}T_g.
\end{equation*}
Its dual (\cite[Theorem 6.4.6]{MR2254662}) with respect to the Petersson inner product is given by
\begin{equation*}
T_n^\ast=\frac{1}{n}\sum_{\substack{g=\text{diag}(y_1,y_2,y_3) \\
y_1\mid y_2\mid y_3 \\
y_1y_2y_3=n}}T_{g^{-1}}.
\end{equation*}
\par
Applying the amplification method requires being able to linearize the composition of several Hecke operators. The different required formulas proved in \cite{MR3525536} are encapsulated in the proposition.
\begin{proposition}[R.~Holowinsky-G.~Ricotta-E.~Royer (\cite{MR3525536})]\label{propo_linear}
Let $p$ and $q$ be two prime numbers.
\begin{eqnarray*}
T_p\circ T_q & = & \frac{1}{pq}T_{\diag(1,1,pq)}+\delta_{p=q}\frac{p+1}{p^2}T_{\diag(1,p,p)}, \\
T_p^\ast\circ T_q & = & \frac{1}{pq}T_{\diag(1,p,pq)}+\delta_{p=q}\frac{p^2+p+1}{p^2}\id, \\
T_p^\ast\circ T_q^\ast & = & \frac{1}{pq}T_{\diag(1,pq,pq)}+\delta_{p=q}\frac{p+1}{p^2}T_{\diag(1,1,p)}. \end{eqnarray*}
\begin{multline*}
T_p\circ\left(T_q\circ T_q^\ast-\id\right)=\frac{q+1}{pq^2}T_{\diag(1,1,p)}+\frac{1}{pq^2}T_{\diag(1,q,pq^2)} \\
+\delta_{p=q}\left(\frac{p+1}{p^3}T_{\diag(1,p^2,p^2)}+\frac{p+1}{p^2}T_{\diag(1,1,p)}\right).
\end{multline*}
\begin{multline*}
T_p^\ast\circ\left(T_q\circ T_q^\ast-\id\right)=\frac{q+1}{pq^2}T_{\diag(1,p,p)}+\frac{1}{pq^2}T_{\diag(1,pq,pq^2)} \\
+\delta_{p=q}\left(\frac{p+1}{p^3}T_{\diag(1,1,p^2)}+\frac{p+1}{p^2}T_{\diag(1,p,p)}\right).
\end{multline*}
\begin{multline*}
\left(T_p\circ T_p^\ast-\id\right)\circ\left(T_q\circ T_q^\ast-\id\right)=\frac{1}{p^2q^2}T_{\diag(1,pq,p^2q^2)}+\frac{q+1}{p^2q^2}T_{\diag(1,p,p^2)} \\
+\frac{p+1}{p^2q^2}T_{\diag(1,q,q^2)}+\frac{(p+1)(q+1)}{p^2q^2}\id \\
+\delta_{p=q}\left(\frac{p+1}{p^4}T_{\diag(1,p^3,p^3)}+\frac{p+1}{p^4}T_{\diag(1,1,p^3)}\right) \\
+\delta_{p=q}\left(\frac{(p+1)(2p-1)}{p^4}T_{\diag(1,p,p^2)}+\frac{p(p+1)(1+p+p^2)}{p^4}\id\right).
\end{multline*}
Moreover,
\begin{eqnarray*}
T_{p,1}=T_{1,p}^\ast & = & T_p, \\
T_{p,1}^\ast=T_{1,p} & = & T_p^\ast, \\
T_{p,p}=T_{p,p}^\ast & = & T_p\circ T_p^\ast-\id.
\end{eqnarray*}
\end{proposition}
Recall that the Hecke algebra is isomorphic to the algebra of double $GL_3(\Z)$-cosets where the multiplication law is defined in \cite{MR1291394}. The previous proposition follows from an explicit computation of the multiplication of the corresponding double cosets.
\subsection{Constructing an amplifier}%
In this section, we will choose the set $I$ and the coefficients $\alpha_{m,n}$, $(m,n)\in I$ occurring in \eqref{eq_amplified_pre_trace}.
\par
Let us construct a relevant $GL(3)$ amplifier, based on the identity
\begin{equation}\label{eq_multi}
a_{j_0}(1,p)a_{j_0}(p,1)-a_{j_0}(p,p)=1
\end{equation}
where $a_{j_0}(m,n)$ stands for the $(m,n)$-th Fourier coefficient of $\Phi_{j_0}$. Let $L\geq 1$ be a parameter, whose value will be determined later on (a positive power of $T$). Let us choose
\begin{equation}\label{eq_I_choice}
I\coloneqq\left\{(p,1), (1,p), (p,p), L\leq p\leq 2L, p \text{ prime}\right\}
\end{equation}
and
\begin{equation}\label{eq_choice_ampli}
\alpha_{m,n}\coloneqq\begin{cases}
a_{j_0}(1,p) & \text{if $L\leq m=p\leq 2L$ is a prime and $n=1$,} \\
a_{j_0}(p,1) & \text{if $m=1$ and $L\leq n=p\leq 2L$ is a prime,} \\
-2 & \text{if $L\leq m=n=p\leq 2L$ are the same prime,} \\
0 & \text{otherwise}
\end{cases}
\end{equation}
such that
\begin{align*}
A_{j_0}(\alpha) & =2\sum_{L\leq p\leq 2L}\left(a_{j_0}(1,p)a_{j_0}(p,1)-a_{j_0}(p,p)\right) \\
& =2\sum_{L\leq p\leq 2L}1
\end{align*}
satisfies
\begin{equation}\label{eq_bound_ampli}
A_{j_0}(\alpha)\gg_\epsilon L^{1-\epsilon}
\end{equation}
by \eqref{eq_multi}. 
\section{Test functions in the pre-trace formula}\label{sec_test}%
\subsection{On the cuspidal spectrum of $X$}
Let $\Phi$ be a Hecke-Maass cusp form of type $\left(\nu_1,\nu_2\right)\in\C^2$. Its archimedean Langlands parameters are
\begin{equation*}
\left(\alpha_1,\alpha_2,\alpha_3\right)=\left(2\nu_1+\nu_2,-\nu_1+\nu_2,-\nu_1-2\nu_2\right)
\end{equation*}
and the element of $\mathfrak{a}_\C^\ast/W$ corresponding to $\Phi$ is
\begin{equation*}
\lambda_\Phi=3\nu_1\lambda_1+3\nu_2\lambda_2.
\end{equation*}
Let us denote by $\Lambda$ the set of these linear forms. The Laplacian eigenvalue of $\Phi$ is
\begin{equation*}
1-3\nu_1^2-3\nu_1\nu_2-3\nu_2^2=1-\frac{1}{2}\left(\alpha_1^2+\alpha_2^2+\alpha_3^2\right).
\end{equation*}
The Jacquet-Shalika bound towards the Ramanujan-Petersson-Selberg conjecture asserts that
\begin{equation*}
\max_{1\leq i\leq 3}\left\vert\Re{\left(\alpha_i\right)}\right\vert\leq\frac{1}{2}
\end{equation*}
and the unitaricity condition tells us that
\begin{equation*}
\left\{\alpha_1,\alpha_2,\alpha_3\right\}=\left\{-\overline{\alpha_1},-\overline{\alpha_2},-\overline{\alpha_3}\right\}.
\end{equation*}
Both previous facts ensure that either
\begin{equation*}
\left(\nu_1,\nu_2\right)\in\left(i\R\right)^2,
\end{equation*}
in which case $\Phi$ is said to be tempered or
\begin{equation*}
\left(\nu_1,\nu_2\right)=\left(\frac{2\sigma}{3},-\frac{\sigma}{3}+it\right)
\end{equation*}
with $\sigma$ and $t$ in $\R$ with $\left\vert\sigma\right\vert\leq 1/2$, in which case $\Phi$ is said to be exceptionnal.
\subsection{Construction of a relevant test function on the spectral side}\label{subsec_construct}%
In this section, we will design the function $\mathcal{H}(k)$ occurring in \eqref{eq_amplified_pre_trace}.
\par
If $F=\{a\in A, d(a,I)\geq 1\}$ then $F$ is a closed subset of $G$, which does not contain $I$. By the properties of the distance function, $g$ in $KFK$ also satisfies $d(g, I) \geq 1$. Thus, one can find a Weyl-invariant symmetric open neighborhood $O$ of $I$ in $G$ and a small enough positive real number $\delta$ satisfying
\begin{equation*}
I\in O\subset A(\delta)=\left\{a\in A,\left\vert\left\vert\log{a}\right\vert\right\vert\leq\delta\right\}\subset G\setminus KFK
\end{equation*}
and $KA(\delta)K\subset G\setminus KFK=\{g\in G, d(g,I)<1\}$.
\par
\begin{figure}[H]
\centering 
{\includegraphics{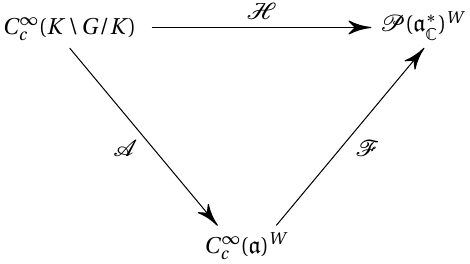}} 
\caption{The Paley-Wiener theorem}
\label{fig:helgason} 
\end{figure}

The Paley-Wiener theorem asserts that the diagram given in figure \ref{fig:helgason} is a commutative diagram of isomorphisms of topological algebras. In this diagram, $\mathcal{H}$ is the Helgason transform, $\mathcal{F}$ is the Fourier transform and $\mathcal{A}$ is the Abel transform. Of course, $C^\infty_c(\mathfrak{a})^W$ can be identified to $C^\infty_c(A)^W$, via the exponential map. R. Gangolli proved a refined version in \cite{MR0289724} of the Paley-Wiener theorem, which says that if $g$ belongs to $C^\infty_c(A(\delta))^W$ then $\mathcal{A}^{-1}(g)$ belongs to $C^\infty_c(KA(\delta)K)\subset C^\infty_c(G\setminus KFK)$.
\par
Both previous paragraphs imply that there exists a Weyl-invariant symmetric open neighborhood $U$ of $0$ in $\mathfrak{a}$ such that
\begin{equation*}
\forall g\in C^\infty_c(U)^W, \quad\mathcal{A}^{-1}(g)\in C^\infty_c(G\setminus KFK)
\end{equation*} 
and $\vert\vert H\vert\vert\leq 1/3$ for $H$ in $U$. 
\par
\begin{figure}[H]
\centering 
\scalebox{0.5} 
{\includegraphics{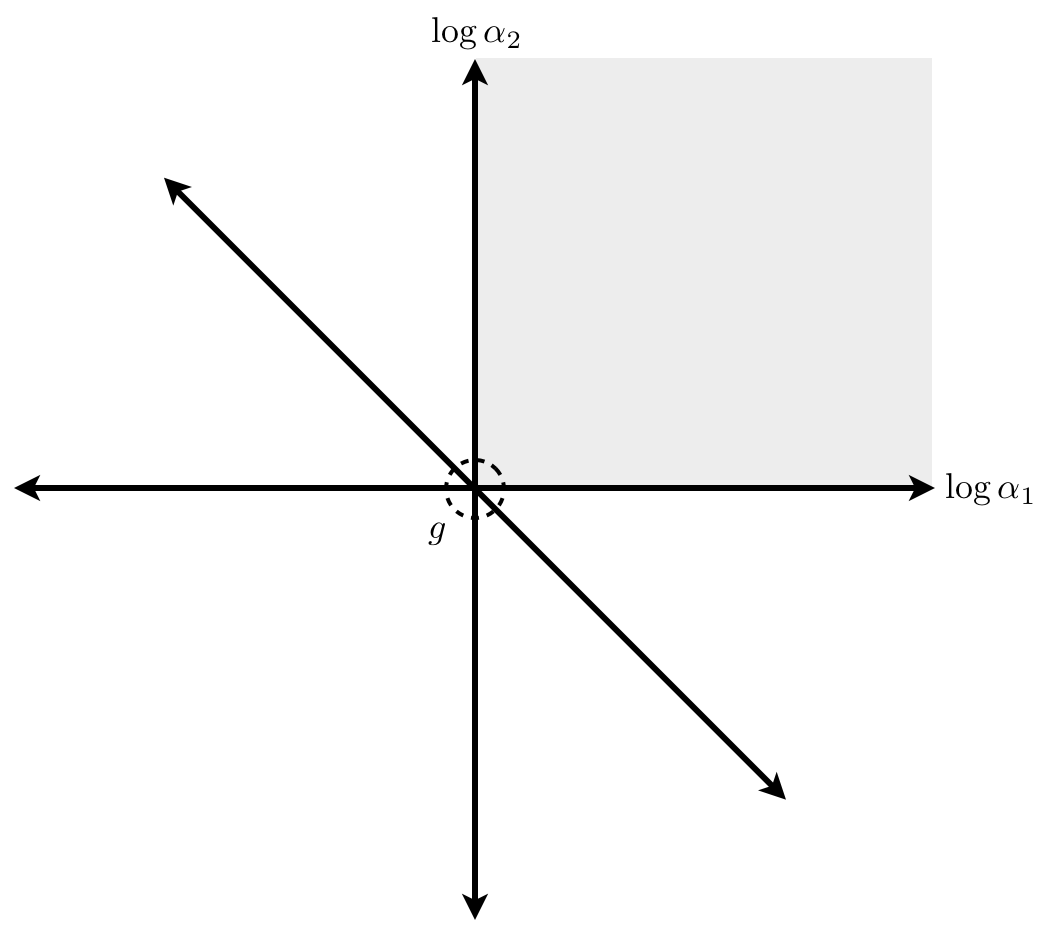}} 
\caption{Test function $g\in C_c^{\infty}(\mathfrak{a})^W$}
\label{fig:gtest} 
\end{figure}

Let us fix $U'$ a Weyl-invariant symmetric open neighborhood of $0$ in $\mathfrak{a}$ satisfying $U'+U'\subset U$. Let us also fix a real non-negative symmetric function $g$ in $C^\infty_c(U')^W$ normalized by $\int_{h\in\mathfrak{a}}g(h)\mathrm{d}H=2$. See figure \ref{fig:gtest}.
By \cite[Lemma 6.2]{MR532745}, the function $\mathcal{F}(g)$ in $\mathcal{P}(\mathfrak{a}_{\mathbb{C}}^\ast)^W$ is even , real-valued\footnote{But not on $\mathfrak{a}_\C^\ast$.} on the spectrum $\Lambda$ of $X$  and satisfies
\begin{equation*}
\forall\lambda\in\mathfrak{a}_{\mathbb{C}}^\ast,\quad\vert\vert\lambda\vert\vert\leq 1\Rightarrow \vert\mathcal{F}(g)(\lambda)\vert\geq 1.
\end{equation*}
Recall that the Paley-Wiener condition means that
\begin{equation*}
\forall\lambda\in\mathfrak{a}_{\mathbb{C}}^\ast, \forall m\geq 0,\quad\left\vert\mathcal{F}(g)(\lambda)\right\vert\leq c_m(g)\frac{\exp{\left(\delta\vert\vert\lambda_{\mathbb{R}}\vert\vert\right)}}{\left(1+\vert\vert\lambda\vert\vert\right)^m}.
\end{equation*}
Briefly speaking, $\mathcal{F}(g)$ is a real bump function over $0$.
\par
In order to restore the positivity, let us define $h=g\ast g$ such 
that $\mathcal{F}(h)=\mathcal{F}(g)^2$. By \cite[Lemma 6.3]{MR532745}, the function $h$ in $C^\infty_c(U)^W$ is real symmetric and its Fourier transform $\mathcal{F}(h)$, which belongs to $\mathcal{P}(\mathfrak{a}_{\mathbb{C}}^\ast)^W$, is a non-negative\footnote{But not on $\mathfrak{a}_\C^\ast$.} function on the spectrum $\Lambda$ of $X$ satisfying
\begin{equation*}
\forall\lambda\in\mathfrak{a}_{\mathbb{C}}^\ast,\quad\vert\vert\lambda\vert\vert\leq 1\Rightarrow \left\vert\mathcal{F}(h)(\lambda)\right\vert\geq 1.
\end{equation*}
The Paley-Wiener condition becomes
\begin{equation*}
\forall\lambda\in\mathfrak{a}_{\mathbb{C}}^\ast, \forall m\geq 0,\quad\mathcal{F}(h)(\lambda)\leq d_m(g)\frac{\exp{\left(2\delta\vert\vert\lambda_{\mathbb{R}}\vert\vert\right)}}{\left(1+\vert\vert\lambda\vert\vert\right)^m}.
\end{equation*}
Thus, $\mathcal{F}(h)$ is a bump function over $0$ non-negative on the spectrum $\Lambda$ of $X$.
\par
We would like to construct a bump function over the spectral parameter of our favorite tempered Hecke-Maass cusp form $\Phi_{0}$. By the previous section, the element $\lambda_{\Phi_0}$ of $\mathfrak{a}_\C^\ast$ associated to $\Phi_0$ is given by
\begin{equation*}
\lambda_{\Phi_0}=3\nu_{0,1}\lambda_1+3\nu_{0,2}\lambda_2\in i\mathfrak{a}^\ast
\end{equation*}
where $(\nu_{0,1},\nu_{0,2})$ is the type of $\Phi_0$, which belongs to $i\R^2$ by the temperedness condition on $\Phi_0$. Let us define
\begin{equation}\label{eq_muT}
\mu_T=3iT\lambda_1+3iT\lambda_2
\end{equation}
and
\begin{equation*}
h_T=e^{-\mu_T}h\rightsquigarrow\mathcal{F}(h_T)(\lambda)=\mathcal{F}(h)(\lambda-\mu_T).
\end{equation*}
This function $h_T$ belongs to $C_c^\infty(U)$ and its Fourier transform satisfies
\begin{equation*}
\forall\lambda\in\mathfrak{a}_{\mathbb{C}}^\ast,\quad\vert\vert\lambda-\mu_T\vert\vert\leq 1\Rightarrow \left\vert\mathcal{F}(h_T)(\lambda)\right\vert\geq 1.
\end{equation*}
The Paley-Wiener condition becomes
\begin{equation}\label{eq_PW_ht}
\forall\lambda\in\mathfrak{a}_{\mathbb{C}}^\ast, \forall m\geq 0,\quad\mathcal{F}(h_T)(\lambda)\leq d_m(g)\frac{\exp{\left(2\delta\vert\vert\rho\vert\vert\right)}}{\left(1+\vert\vert\lambda-\mu_T\vert\vert\right)^m}.
\end{equation}
This follows from the Paley-Wiener condition for $h$ and the fact that $(\lambda-\mu_T)_{\mathbb{R}}=\lambda_{\mathbb{R}}$ with $\vert\vert\lambda_{\mathbb{R}}\vert\vert\leq\vert\vert\rho\vert\vert$ by \cite[Proposition 3.4]{MR532745}. Thus, $\mathcal{F}(h_T)$ is a bump function over $\mu_T$ non-negative on the spectrum $\Lambda$ of $X$.
\par
With $h_T$ not Weyl-invariant, it seems natural to define
\begin{equation*}
h_T^W(H)=\sum_{w\in W}h_T(w.h)=h(H)\sum_{w\in W}e^{-\mu_T(w.h)}
\end{equation*}
whose Fourier transform is given by
\begin{equation*}
\mathcal{F}\left(h_T^W\right)(\lambda)=\sum_{w\in W}\mathcal{F}(h)(\lambda-w.\mu_T).
\end{equation*}
The previous paragraphs imply that $h_T^W$ belongs to $C_c^\infty(U)^W$. In 
particular, $\mathcal{A}^{-1}(h_T^W)$ is supported in the compact 
set $G\setminus KFK$ which does not depend on $T$. The Fourier transform of
$h_T^W$ is non-negative on the spectrum $\Lambda$ of $X$ and satisfies for 
$\lambda\in\mathfrak{a}_{\mathbb{C}}^\ast$
\begin{equation*}
\left\vert\mathcal{F}\left(h_T^W\right)(\lambda)\right\vert\geq 1.
\end{equation*}
as soon as there exists $w$ in $W$ with $\vert\vert\lambda-w.\mu_T\vert\vert\leq 1$.
\par
\begin{figure}[H]
\centering 
\scalebox{0.5} 
{\includegraphics{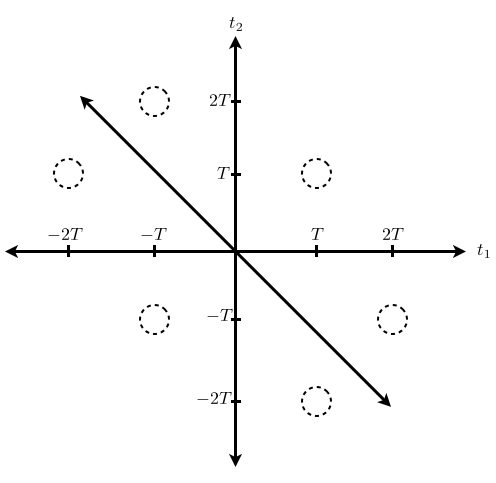}} 
\caption{Test function $\mathcal{F}\left(h_T^W\right)\in \mathcal{P}(\mathfrak{a}_\mathbb{C}^{\ast})^W$}
\label{fig:htest} 
\end{figure}

This function $\mathcal{F}(h_T^W)$ is the 
Weyl-invariant bump function non-negative on the spectrum $\Lambda$ of $X$ we were looking at (see figure \ref{fig:htest}). 
In other words, $\mathcal{H}(k)=\mathcal{F}(h_T^W)$ in 
\eqref{eq_amplified_pre_trace}, and $k = \mathcal{A}^{-1}(h_T^W)$.
\subsection{Estimates for the inverse Helgason transform of our test function}%
The \emph{spherical function} of parameter $s\in\C^2$ is defined by
\[
  \varphi_s(g) = \int_{k\in K}(p_s\delta^{1/2})(\IwA(kg))\dd k
\]
for $g$ in $G$ with the Haar measure on $K$ normalized so that $K$ has measure 
one. The spherical function $\varphi_s$ is a bi-$K$-invariant function on $G$, Weyl-invariant
in its parameter $s$ and satisfies
$\varphi_s(I) = 1$. We will also write $\varphi_\lambda$ where the association between
$\lambda$ and $s$ is as in \eqref{eq_lambda_s}. The oscillatory integral which forms the spherical
function has been studied by
many authors, including J.~J.~Duistermaat, J.~A.~C.~Kolk and V.~S.~Varadarajan \cite{MR532745}, 
V.~Blomer and A.~Pohl \cite{BlPo} and S.~Marshall \cite{Ma3}. We will rely
on the result of S.~Marshall, which we restate below just for $GL(3)$ and in
our notation, though his
result is for semisimple and noncompact groups with finite center. Define
the \emph{singular set} in $i\mathfrak{a}^\ast$ to be 
\[
  \left\{\lambda\in i\mathfrak{a}^\ast, B(\alpha_j^+, \lambda)\in \pi i\Z\textnormal{ for some } j=1,2,3\right\}.
\]
\begin{proposition}[S.~Marshall, Theorem 1.3, \cite{Ma3}]\label{propo_marshall}
  Let $B\subset A$ be a compact set and let $B^\ast\subset i\mathfrak{a}^\ast$
  be a compact set which does not intersect the singular set. Then
  \begin{equation}
    \varphi_{\exp(T\lambda)}(a) \ll_{B,B^\ast} \prod_{j=1}^3 (1 + T|\log \alpha_j(a)|)^{-1/2}
    \label{eq_marshall_bound}
  \end{equation}
  for any $a$ in $B$ and $\lambda$ in $B^\ast$.
\end{proposition}
\par
\begin{figure}[H]
\centering 
\scalebox{0.5} 
{\includegraphics{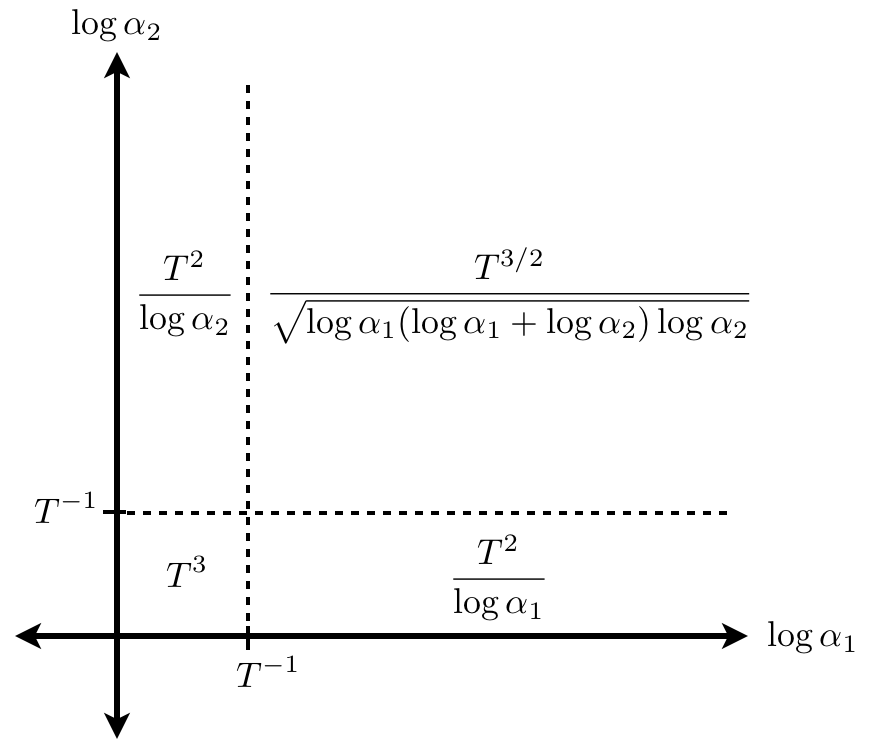}} 
\caption{Bounds for the inverse Helgason transform $k=\mathcal{H}^{-1}\left(\mathcal{F}\left(h_T^W\right)\right)$ (up to $T^\epsilon$).}
\label{fig:kbounds} 
\end{figure}

The inverse Helgason transform, also called the inverse spherical transform, 
is given by 
\begin{equation}\label{eq_inversion_1}
k(a)=\mathcal{H}^{-1}\left(\mathcal{F}\left(h_T^W\right)\right)(a)=\int_{t\in\R^2}\mathcal{F}\left(h_T^W\right)(t)\varphi_{it}(a)\frac{\dd t}{\abs{c_3(t)}^2},
\end{equation}
the measure being the Plancherel one, where $c_3$ stands for the 
Harish-Chandra $c$-function. The required estimates for the inverse Helgason 
transform $k$ of our test function $\mathcal{F}(h_T^W)$ constructed in 
Subsection \ref{subsec_construct} which will enable us to estimate the geometric side of the 
amplified pre-trace formula \eqref{eq_amplified_pre_trace} are given in the following 
proposition.
\begin{proposition}\label{propo_estimates}
Let $a$ be an element in a compact subset of $A$.
\begin{itemize}
\item
If $a$ belongs to the closure of the positive Weyl chamber $A_+$ then
\begin{equation*}
  k(a)=\mathcal{H}^{-1}\left(\mathcal{F}\left(h_T^W\right)\right)(a)\ll_{\epsilon} T^{3+\epsilon}.
\end{equation*}
\item
If $a$ belongs to the positive Weyl chamber $A_+$ then
\begin{equation*}
  k(a)=\mathcal{H}^{-1}\left(\mathcal{F}\left(h_T^W\right)\right)(a)\ll_{\epsilon}\frac{T^{3/2+\epsilon}}{\sqrt{\left(\alpha_1(a)^2-1\right)\left(\alpha_2(a)^2-1\right)\left(\alpha_3(a)^2-1\right)}}.
\end{equation*}
\item
If $a$ satisfies $1\leq\alpha_1(a)\leq 1+O(1)/T$ and $\alpha_2(a)\geq 1+O(1)/T$ then
\begin{equation*}
  k(a)=\mathcal{H}^{-1}\left(\mathcal{F}\left(h_T^W\right)\right)(a)\ll_{\epsilon}\frac{T^{2+\epsilon}}{\alpha_2(a)^2-1}.
\end{equation*}
\item
If $a$ satisfies $\alpha_1(a)\geq 1+O(1)/T$ and $1\leq\alpha_2(a)\leq 1+O(1)/T$ then
\begin{equation*}
  \mathcal{H}^{-1}\left(\mathcal{F}\left(h_T^W\right)\right)(a)\ll_{\epsilon}\frac{T^{2+\epsilon}}{\alpha_1(a)^2-1}.
\end{equation*}
\end{itemize}

\end{proposition}
Altogether, the bounds given in this proposition are summarized in the figure \ref{fig:kbounds}.
\begin{proof}
  By the Weyl-invariance of both $\varphi_{it}$ in its parameter and of the Plancherel measure, and by the construction of the test function $\mathcal{F}(h_T^W)$,
  \[
    k(a) = 6\int_{t\in\R^2}\mathcal{F}\left(h_T\right)(t)
    \varphi_{it}(a)\frac{\dd t}{|c_3(it)|^2}.
  \]
  The spherical function satisfies $|\varphi_{it}(a)|$ is bounded while
  the Harish-Chandra $c$-function satisfies (see \cite[Chapter 5, Theorem 6.4]{MR2166237})
  \[
    |c_3(i(t_1, t_2))|^{-2} = \frac{\pi}{12}t_1t_2(t_1 + t_2)\tanh\left(\frac{\pi}{2} t_1\right)\tanh\left(\frac{\pi}{2} t_2\right)
    \tanh\left(\frac{\pi}{2}(t_1 + t_2)\right),
  \]
  and thus grows polynomially in $t$.
  The Paley-Wiener estimate \eqref{eq_PW_ht} of arbitrary polynomial decay 
  for the test function away from $\mu_T$ implies that for any positive 
  integer $m$,
  \[
    k(a) = 6\int_{B_{T^\epsilon}(\mu_T)}
    \mathcal{F}\left(h_T\right)(t)\varphi_{it}(a)
    \frac{\dd t}{|c_3(it)|^2} + O_m(T^{-m})
  \]
	where $\mu_T$ is defined in \eqref{eq_muT} and $B_{T^\epsilon}(\mu_T)$ stands for a ball of center $(T,T)$ and radius $T^\epsilon$.
  In $B_{T^\epsilon}(\mu_T)$, $|c_3(it)|^2 \ll T^{3}$ and 
  $\mathcal{F}\left(h_T\right)$ is
  bounded. We will now see that the conditions for S.~Marshall's bound 
  \eqref{eq_marshall_bound} are met. 
  Let $B$ be a compact set which contains the support of $k$ for all $T$. 
  Such sets exist by the construction of $k$. Since
  $i(\lambda_1 + \lambda_2)$ is not in the singular set, it is possible to
  take $B^\ast$ to be a closed ball around $i(\lambda_1 + \lambda_2)$ that is
  disjoint from the singular set. Then for $T$ sufficiently large, $TB^\ast$ 
  will contain a ball of radius $T^\epsilon$ around $\mu_T$. With these 
  choices of the sets
  $B$ and $B^\ast$ made, \eqref{eq_marshall_bound} applies.
  Taylor expanding $\log\alpha_i$ at $\alpha_i = 1$ if
  $a$ is near a wall of the Weyl chamber gives the denominators in the 
  proposition.
\end{proof}
\section{First estimate for the geometric side of the amplified pre-trace formula}\label{sec_first_estimate}%
This section is devoted to the proof of the following first estimate for $\left\vert\Phi_{j_0}(z)\right\vert$. Let us define
\begin{equation}\label{eq_Kmn}
K_{\ell,n}(z)\coloneqq\sum_{\rho\in GL_3(\Z)\diag(1,\ell,n)GL_3(\Z)/\left\{\pm 1\right\}}\left\vert k\left(\frac{1}{\text{det}(\rho)^{1/3}}z^{-1}\rho z\right)\right\vert
\end{equation}
for any positive integer $n$, any positive integer $\ell$ dividing $n$ and any $z$ in $X$ and where
\begin{equation*}
k=\mathcal{H}^{-1}\left(\mathcal{F}\left(h_T^W\right)\right).
\end{equation*}
\begin{proposition}\label{propo_first_11}
Let $z$ be in $X$. One has
\begin{multline}\label{eq_geo_side_1}
L^{2-\epsilon}\left\vert\Phi_{j_0}(z)\right\vert^2\ll_\epsilon\sum_{L\leq p,q\leq 2L}\frac{\left\vert\alpha_{p,1}\alpha_{q,1}\right\vert}{pq}K_{q,pq}(z)+\sum_{L\leq p\leq 2L}\frac{\left\vert\alpha_{p,1}\right\vert^2(p^2+p+1)}{p^2}K_{1,1}(z) \\
+\sum_{L\leq p,q\leq 2L}\frac{\left\vert\alpha_{p,1}\alpha_{1,q}\right\vert}{pq}K_{1,pq}(z)+\sum_{L\leq p\leq 2L}\frac{\left\vert\alpha_{p,1}\right\vert^2(p+1)}{p^2}K_{p,p}(z) \\
+\sum_{L\leq p,q\leq 2L}\frac{\left\vert\alpha_{p,1}\right\vert(q+1)}{pq^2}K_{1,p}(z)+\sum_{L\leq p,q\leq 2L}\frac{\left\vert\alpha_{p,1}\right\vert}{pq^2}K_{q,pq^2}(z) \\
+\sum_{L\leq p\leq 2L}\frac{\left\vert\alpha_{p,1}\right\vert(p+1)}{p^3}K_{p^2,p^2}(z)+\sum_{L\leq p\leq 2L}\frac{\left\vert\alpha_{p,1}\right\vert(p+1)}{p^2}K_{1,p}(z) \\
+\sum_{L\leq p,q\leq 2L}\frac{1}{p^2q^2}K_{pq,p^2q^2}(z)+\sum_{L\leq p,q\leq 2L}\frac{q+1}{p^2q^2}K_{p,p^2}(z) \\
+\sum_{L\leq p,q\leq 2L}\frac{p+1}{p^2q^2}K_{q,q^2}(z)+\sum_{L\leq p,q\leq 2L}\frac{(p+1)(q+1)}{p^2q^2}K_{1,1}(z) \\
+\sum_{L\leq p\leq 2L}\frac{p+1}{p^4}K_{p^3,p^3}(z)+\sum_{L\leq p\leq 2L}\frac{p+1}{p^4}K_{1,p^3}(z) \\
+\sum_{L\leq p\leq 2L}\frac{(p+1)(2p-1)}{p^4}K_{p,p^2}(z)+\sum_{L\leq p\leq 2L}\frac{p(p+1)(1+p+p^2)}{p^4}K_{1,1}(z)
\end{multline}
where all the summations are over prime numbers.
\end{proposition}
The quantities $K_{\ell,n}(z)$ will be bounded thanks to Proposition \ref{propo_estimates} and a counting lemma given in the next section.
\begin{proof}[\proofname{} of Proposition \ref{propo_first_11}]
The amplifier defined in \eqref{eq_choice_ampli} satisfies 
\begin{equation}\label{eq_ampli_1}
\overline{\alpha_{m,n}}=\alpha_{n,m}
\end{equation}
and
\begin{equation}\label{eq_ampli_2}
\overline{\alpha_{m,m}}=\alpha_{m,m}
\end{equation}
for any $(m,n)\in I$, the set defined in \eqref{eq_I_choice}.
\par
Let $z'$ be in $X$ and define
\begin{equation*}
S\coloneqq\sum_{j\geq 0}\left\vert\sum_{(m,n)\in I}\alpha_{m,n}a_j(m,n)\right\vert^2\hat{h}(\nu_j)\Phi_j(z)\overline{\Phi_j}(z^\prime).
\end{equation*}
Expanding the square,
\begin{equation*}
S=\sum_{k=1}^9S_k(g,g^\prime)
\end{equation*}
where
\begin{eqnarray*}
S_1(z,z^\prime) & = & \sum_{p,q\sim L}\alpha_{p,1}\overline{\alpha_{q,1}}\sum_{j\geq 0}\hat{h}(\nu_j)a_j(p,1)\overline{a_j(q,1)}\Phi_j(z)\overline{\Phi_j}(z^\prime), \\
S_2(z,z^\prime) & = & \sum_{p,q\sim L}\alpha_{p,1}\overline{\alpha_{1,q}}\sum_{j\geq 0}\hat{h}(\nu_j)a_j(p,1)\overline{a_j(1,q)}\Phi_j(z)\overline{\Phi_j}(z^\prime), \\
S_3(z,z^\prime) & = & \sum_{p,q\sim L}\alpha_{p,1}\overline{\alpha_{q,q}}\sum_{j\geq 0}\hat{h}(\nu_j)a_j(p,1)\overline{a_j(q,q)}\Phi_j(z)\overline{\Phi_j}(z^\prime)
\end{eqnarray*}
and
\begin{eqnarray*}
S_4(z,z^\prime) & = & \sum_{p,q\sim L}\alpha_{1,p}\overline{\alpha_{q,1}}\sum_{j\geq 0}\hat{h}(\nu_j)a_j(1,p)\overline{a_j(q,1)}\Phi_j(z)\overline{\Phi_j}(z^\prime), \\
S_5(z,z^\prime) & = & \sum_{p,q\sim L}\alpha_{1,p}\overline{\alpha_{1,q}}\sum_{j\geq 0}\hat{h}(\nu_j)a_j(1,p)\overline{a_j(1,q)}\Phi_j(z)\overline{\Phi_j}(z^\prime), \\
S_6(z,z^\prime) & = & \sum_{p,q\sim L}\alpha_{1,p}\overline{\alpha_{q,q}}\sum_{j\geq 0}\hat{h}(\nu_j)a_j(1,p)\overline{a_j(q,q)}\Phi_j(z)\overline{\Phi_j}(z^\prime)
\end{eqnarray*}
and
\begin{eqnarray*}
S_7(z,z^\prime) & = & \sum_{p,q\sim L}\alpha_{p,p}\overline{\alpha_{q,1}}\sum_{j\geq 0}\hat{h}(\nu_j)a_j(p,p)\overline{a_j(q,1)}\Phi_j(z)\overline{\Phi_j}(z^\prime), \\
S_8(z,z^\prime) & = & \sum_{p,q\sim L}\alpha_{p,p}\overline{\alpha_{1,q}}\sum_{j\geq 0}\hat{h}(\nu_j)a_j(p,p)\overline{a_j(1,q)}\Phi_j(z)\overline{\Phi_j}(z^\prime), \\
S_9(z,z^\prime) & = & \sum_{p,q\sim L}\alpha_{p,p}\overline{\alpha_{q,q}}\sum_{j\geq 0}\hat{h}(\nu_j)a_j(p,p)\overline{a_j(q,q)}\Phi_j(z)\overline{\Phi_j}(z^\prime).
\end{eqnarray*}
One can check that
\begin{eqnarray*}
\overline{S_9(z,z^\prime)} & = & S_9(z^\prime,z) \\
\overline{S_1(z,z^\prime)} & = & S_5(z^\prime,z) \\
\overline{S_2(z,z^\prime)} & = & S_4(z^\prime,z) \\
\overline{S_3(z,z^\prime)} & = & S_6(z^\prime,z) \\
\overline{S_7(z,z^\prime)} & = & S_8(z^\prime,z)
\end{eqnarray*}
by \eqref{eq_ampli_1} and \eqref{eq_ampli_2}. Thus,
\begin{equation*}
S=\sum_{k=1}^4\left(T_k(z,z^\prime)+\overline{T_k(z^\prime,z)}\right)+T_5(z,z^\prime)
\end{equation*}
where
\begin{eqnarray*}
T_1(z,z^\prime) & = & \sum_{p,q\sim L}\alpha_{p,1}\overline{\alpha_{q,1}}\sum_{j\geq 0}\hat{h}(\nu_j)a_j(p,1)\overline{a_j(q,1)}\Phi_j(z)\overline{\Phi_j}(z^\prime), \\
T_2(z,z^\prime) & = & \sum_{p,q\sim L}\alpha_{p,1}\overline{\alpha_{1,q}}\sum_{j\geq 0}\hat{h}(\nu_j)a_j(p,1)\overline{a_j(1,q)}\Phi_j(z)\overline{\Phi_j}(z^\prime), \\
T_3(z,z^\prime) & = & -2\sum_{p,q\sim L}\alpha_{p,1}\sum_{j\geq 0}\hat{h}(\nu_j)a_j(p,1)a_j(q,q)\Phi_j(z)\overline{\Phi_j}(z^\prime), \\
T_4(z,z^\prime) & = & -2\sum_{p,q\sim L}\overline{\alpha_{q,1}}\sum_{j\geq 0}\hat{h}(\nu_j)a_j(p,p)\overline{a_j(q,1)}\Phi_j(z)\overline{\Phi_j}(z^\prime), \\
T_5(z,z^\prime) & = & 4\sum_{p,q\sim L}\sum_{j\geq 0}\hat{h}(\nu_j)a_j(p,p)a_j(q,q)\Phi_j(z)\overline{\Phi_j}(z^\prime).
\end{eqnarray*}
One can check that
\begin{equation*}
\overline{T_3(z,z^\prime)}=T_4(z^\prime,z)
\end{equation*}
such that
\begin{equation*}
S=\sum_{k=1}^2\left(U_k(z,z^\prime)+\overline{U_k(z^\prime,z)}\right)+2\left(U_3(z,z^\prime)+\overline{U_3(z^\prime,z)}\right)+U_4(z,z^\prime)
\end{equation*}
where
\begin{eqnarray*}
U_1(z,z^\prime) & = & \sum_{p,q\sim L}\alpha_{p,1}\overline{\alpha_{q,1}}\sum_{j\geq 0}\hat{h}(\nu_j)a_j(p,1)\overline{a_j(q,1)}\Phi_j(z)\overline{\Phi_j}(z^\prime) \\
U_2(z,z^\prime) & = & \sum_{p,q\sim L}\alpha_{p,1}\overline{\alpha_{1,q}}\sum_{j\geq 0}\hat{h}(\nu_j)a_j(p,1)\overline{a_j(1,q)}\Phi_j(z)\overline{\Phi_j}(z^\prime) \\
U_3(z,z^\prime) & = & -2\sum_{p,q\sim L}\alpha_{p,1}\sum_{j\geq 0}\hat{h}(\nu_j)a_j(p,1)a_j(q,q)\Phi_j(z)\overline{\Phi_j}(z^\prime) \\
U_4(z,z^\prime) & = & 4\sum_{p,q\sim L}\sum_{j\geq 0}\hat{h}(\nu_j)a_j(p,p)a_j(q,q)\Phi_j(z)\overline{\Phi_j}(z^\prime).
\end{eqnarray*}
Let us define
\begin{equation*}
\varphi(z)=\sum_{j\geq 0}\hat{h}(\nu_j)\Phi_j(z)\overline{\Phi_j}(z^\prime)=\sum_{\gamma\in GL_3(\mathbb{Z})/\{\pm I\}}k(z^{-1}\gamma z^\prime).
\end{equation*}
Now,
\begin{eqnarray*}
U_1(z,z^\prime) & = & \sum_{p,q\sim L}\alpha_{p,1}\overline{\alpha_{q,1}}\left(T_p\circ T_q^\ast\right)\left(\varphi\right)(z) \\
U_2(z,z^\prime) & = & \sum_{p,q\sim L}\alpha_{p,1}\overline{\alpha_{1,q}}\left(T_p\circ T_q\right)\left(\varphi\right)(z) \\
U_3(z,z^\prime) & = & -2\sum_{p,q\sim L}\alpha_{p,1}\left(T_p\circ \left(T_q\circ T_q^\ast-Id\right)\right)\left(\varphi\right)(z) \\
U_4(z,z^\prime) & = & 4\sum_{p,q\sim L}\left(\left(T_q\circ T_q^\ast-Id\right)\circ\left(T_q\circ T_q^\ast-Id\right)\right)\left(\varphi\right)(z).
\end{eqnarray*}
Let us define
\begin{equation*}
K_{m,n}(z,z^\prime)=\sum_{\rho\in GL_3(\Z)\text{diag}(1,m,n)GL_3(\Z)/\{\pm I\}}k\left(\frac{1}{\text{det}(\rho)^{1/3}}z'^{-1}\rho z\right).
\end{equation*}
By the second equation in Proposition \ref{propo_linear} and by \eqref{eq_auto_action},
\begin{equation*}
U_1(z,z^\prime)=\sum_{p,q\sim L}\frac{\alpha_{p,1}\overline{\alpha_{q,1}}}{pq}K_{q,pq}(z,z^\prime)+\sum_{p\sim L}\frac{\abs{\alpha_{p,1}}^2(p^2+p+1)}{p^2}K_{1,1}(z,z^\prime).
\end{equation*}
By the first equation in Proposition \ref{propo_linear} and by \eqref{eq_auto_action},
\begin{equation*}
U_2(z,z^\prime)=\sum_{p,q\sim L}\frac{\alpha_{p,1}\overline{\alpha_{1,q}}}{pq}K_{1,pq}(z,z^\prime)+\sum_{p\sim L}\frac{\alpha_{p,1}^2(p+1)}{p^2}K_{p,p}(z,z^\prime).
\end{equation*}
By the fourth equation in Proposition \ref{propo_linear} and by \eqref{eq_auto_action},
\begin{multline*}
U_3(z,z^\prime)=-2\sum_{p,q\sim L}\frac{\alpha_{p,1}(q+1)}{pq^2}K_{1,p}(z,z^\prime)-2\sum_{p,q\sim L}\frac{\alpha_{p,1}}{pq^2}K_{q,pq^2}(z,z^\prime) \\
-2\sum_{p\sim L}\frac{\alpha_{p,1}(p+1)}{p^3}K_{p^2,p^2}(z,z^\prime)-2\sum_{p\sim L}\frac{\alpha_{p,1}(p+1)}{p^2}K_{1,p}(z,z^\prime).
\end{multline*}
By the sixth equation in Proposition \ref{propo_linear} and by \eqref{eq_auto_action},
\begin{multline*}
U_4(z,z^\prime)=4\sum_{p,q\sim L}\frac{1}{p^2q^2}K_{pq,p^2q^2}(z,z^\prime)+4\sum_{p,q\sim L}\frac{q+1}{p^2q^2}K_{p,p^2}(z,z^\prime) \\
+4\sum_{p,q\sim L}\frac{p+1}{p^2q^2}K_{q,q^2}(z,z^\prime)+4\sum_{p,q\sim L}\frac{(p+1)(q+1)}{p^2q^2}K_{1,1}(z,z^\prime) \\
+4\sum_{p\sim L}\frac{p+1}{p^4}K_{p^3,p^3}(z,z^\prime)+4\sum_{p\sim L}\frac{p+1}{p^4}K_{1,p^3}(z,z^\prime) \\
+4\sum_{p\sim L}\frac{(p+1)(2p-1)}{p^4}K_{p,p^2}(z,z^\prime)+4\sum_{p\sim L}\frac{p(p+1)(1+p+p^2)}{p^4}K_{1,1}(z,z^\prime).
\end{multline*}
Finally, we choose $z^\prime=z$.
\par
The properties of the function $h_T^W$ constructed in the previous section and \eqref{eq_bound_ampli} conclude the proof of this proposition by positivity.
\end{proof}
\section{The counting Lemma}\label{sec_counting}%
\subsection{Preliminary steps}
In this section, $z$ will be in a compact set of $X$, which means that
\begin{equation*}
z=naK=\begin{pmatrix}
1 & x_1 & x_3 \\
& 1 & x_2 \\
& & 1
\end{pmatrix}\begin{pmatrix}
a_1 & & \\
& a_2 & \\
& & a_3
\end{pmatrix}K
\end{equation*}
where
\begin{equation}\label{eq_z_compact}
1\ll x_1, x_2, x_3\ll 1,\quad 1\gg\beta_1\coloneqq\frac{a_1}{a_2}, \beta_2\coloneqq\frac{a_2}{a_3}\geq\frac{\sqrt{3}}{2}.
\end{equation}
\par
In this section, $\rho$ will be an invertible matrix of size $3$, whose Cartan decomposition of $z^{-1}\rho z$ can be written as
\begin{equation*}
z^{-1}\rho z=k_1bk_2=k_1\begin{pmatrix}
b_1 & & \\
& b_2 & \\
& & b_3
\end{pmatrix}k_2 \in KA_+K.
\end{equation*}
By a slight abuse of notations, let us set
\begin{equation*}
\alpha_1=\alpha_1\left(z^{-1}\rho z\right)\coloneqq\alpha_1(b),\quad\alpha_2=\alpha_2\left(z^{-1}\rho z\right)\coloneqq\alpha_2(b)
\end{equation*}
and note that
\begin{equation*}
b_1=\left(n\alpha_1^2\alpha_2\right)^{1/3},\quad b_2=\left(n\frac{\alpha_2}{\alpha_1}\right)^{1/3},\quad b_3=\left(\frac{n}{\alpha_1\alpha_2^2}\right)^{1/3}.
\end{equation*}
\par
Let $M_{\ell,n}(z;\delta_1,\delta_2)$ be the number of matrices
\begin{equation*}
\rho=\begin{pmatrix}
a & b & c \\
d & e & f \\
g & h & j
\end{pmatrix}
\end{equation*}
with integer coefficients satisfying
\begin{equation}\label{eq_constraints}
d(\rho)=\left(d_1(\rho),d_2(\rho),d_3(\rho)\right)=(1,\ell,n),\quad\forall j\in\{1,2\},1\leq\alpha_j\leq 1+\delta_j
\end{equation}
where $\ell$ and $n$ are positive integers with $\ell\mid n$ and $0\leq\delta_1, \delta_2\ll 1$. This section is devoted to the proof of the following proposition. 
\begin{proposition}\label{propo_counting}
Let $z$ be in a compact set of $X$, $0\leq\delta_1, \delta_2\ll 1$ and $\Delta=\delta_1^2+\delta_2^2+\delta_1\delta_2$. One has
\begin{equation*}
M_{\ell,n}(z;\delta_1,\delta_2)\ll_\epsilon n^{1/3+\epsilon}\sum_{\lambda\mid\ell}\frac{1}{\lambda}\left(1+n^{2/3}\left(\sqrt{\Delta}+\Delta\right)^{1/5}\right)^2\left(1+\frac{n^{2/3}\left(\sqrt{\Delta}+\Delta\right)^{1/5}}{\ell/\lambda}\right)\left(1+\frac{n^{1/3}\left(\sqrt{\Delta}+\Delta\right)}{\ell/\lambda}\right)
\end{equation*}
for any $\epsilon>0$.
\end{proposition}
\begin{remark}
The referee kindly pointed us that $\Delta\asymp \delta_1^2+\delta_2^2$ and $\Delta^{1/2}+\Delta\asymp\Delta^{1/2}$ since $0\leq\delta_1, \delta_2\ll 1$. Nevertheless, on the one hand, the statement given in the previous proposition reminds the reader with the distance function given in \eqref{eq_dist} and on the other hand reveals the structure of the proof of this proposition.
\end{remark}
This counting lemma is optimal in the following sense. If $z=I$, the identity matrix, then the number of matrices $\rho$ is bounded by $n^{1/3+\epsilon}$ if $n$ is a cube, which matches the order of magnitude for the number of automorphs of $I$, namely the number of matrices $\rho$ satisfying $\rho K=K$.
\par
The main ingredient in the proof consists in counting integer solutions to equations involving explicit positive definite quadratic forms with real coefficients, which depend on $x_1$, $x_2$, $x_3$ and on the multiplicative roots $\beta_1$ and $\beta_2$. The discriminants of these quadratic forms will be either $\beta_1\geq\sqrt{3}/2>0$ or $\beta_2\geq\sqrt{3}/2>0$, which enables us to approximate them by positive definite quadratic forms with rational coefficients. This Diophantine approximation preliminary step lies at the heart of the proof of the counting lemma proved by V.~Blomer and A.~Pohl (\cite{BlPo}).
\par
Let us fix for now $\rho$, one of these matrices.
\par
One can check that
\begin{equation*}
z^{-1}\rho z=\begin{pmatrix}
a' & b' & c' \\
d' & e' & f' \\
g' & h' & j'
\end{pmatrix}
\end{equation*}
where
\begin{eqnarray*}
a' & = & a-x_1d+xg, \\
b' & = & \frac{(a-x_1d+xg)x_1+b-x_1e+xh}{\beta_1}, \\
c' & = & \frac{(a-x_1d+xg)x_3+(b-x_1e+xh)x_2+c-x_1f+xj}{\beta_3}, \\
d' & = & \beta_1(d-x_2g), \\
e' & = & (d-x_2g)x_1+e-x_2h, \\
f' & = & \frac{(d-x_2g)x_3+(e-x_2h)x_2+f-x_2j}{\beta_2}, \\
g' & = &  \beta_3g, \\
h' & = & \beta_2(gx_1+h), \\
j' & = & gx_3+hx_2+j
\end{eqnarray*}
where $x\coloneqq x_1x_2-x_3$ and $\beta_3\coloneqq\beta_1\beta_2$.
\par
Let us set
\begin{eqnarray*}
\alpha_2 & \coloneqq & dj-fg  \\
\alpha_3 & \coloneqq & dh-eg  \\
\alpha_5 & \coloneqq & aj-cg  \\
\alpha_6 & \coloneqq & ah-bg  \\
\alpha_9 & \coloneqq & ae-bd.
\end{eqnarray*}
\par
The matrix $z^{-1}\rho z$ being close to $n^{1/3}k_1k_2$, let us compute the Frobenius norm of
\begin{equation}\label{eq_Cartan_2}
z^{-1}\rho z-n^{1/3}k_1k_2\eqqcolon\begin{pmatrix}
A & B & C \\
D & E & F \\
G & H & J
\end{pmatrix}.
\end{equation}
By the bi-invariance of the Fronenius norm by orthogonal matrices, one has
\begin{align*}
\abs{\abs{z^{-1}\rho z-n^{1/3}k_1k_2}}_F  & =\sqrt{(b_1-n^{1/3})^2+(b_2-n^{1/3})^2+(b_3-n^{1/3})^2} \\
& \ll n^{1/3}\sqrt{\Delta}
\end{align*}
by \eqref{eq_constraints} and where $\Delta=:\delta_1^2+\delta_2^2+\delta_1\delta_2$. In particular,
\begin{equation}\label{eq_coeff_error}
\abs{A},\dots,\abs{J}\ll n^{1/3}\sqrt{\Delta}
\end{equation}
such that
\begin{equation}\label{eq_coeff'}
\abs{a'},\dots,\abs{j'}\ll n^{1/3}(1+\sqrt{\Delta})\ll n^{1/3}
\end{equation}
since the coefficients of the orthogonal matrix $k_1k_2$ are bounded and
\begin{equation}\label{eq_coeff}
\abs{a},\dots,\abs{j}\ll n^{1/3}
\end{equation}
by the explicit formulas for the coefficients of $z^{-1}\rho z$ and \eqref{eq_z_compact}.
\par
The matrix
\begin{equation*}
k_1k_2=\frac{1}{n^{1/3}}\begin{pmatrix}
a'-A & b'-B & c'-C \\
d'-D & e'-E & f'-F \\
g'-G & h'-H & j'-J
\end{pmatrix}
\end{equation*}
being orthogonal, its rows and columns are orthonormal, which implies
\begin{equation}\label{eq_column_1}
a'^2+d'^2+g'^2=n^{2/3}+O\left(n^{2/3}\left(\sqrt{\Delta}+\Delta\right)\right),
\end{equation}
\begin{equation}\label{eq_column_2}
g'^2+h'^2+j'^2=n^{2/3}+O\left(n^{2/3}\left(\sqrt{\Delta}+\Delta\right)\right),
\end{equation}
and
\begin{equation}\label{eq_column_3}
d'^2+e'^2+f'^2=n^{2/3}+O\left(n^{2/3}\left(\sqrt{\Delta}+\Delta\right)\right)
\end{equation}
by \eqref{eq_coeff_error} and \eqref{eq_coeff'}. In addition, $k_1k_2$ is equal to its comatrix, which implies
\begin{equation}\label{eq_comatrix}
\alpha_2'\coloneqq d'j'-f'g'=\beta_1(\alpha_2+x_2\alpha_3)=-n^{1/3}b'+O\left(n^{2/3}\left(\sqrt{\Delta}+\Delta\right)\right)
\end{equation}
by \eqref{eq_coeff_error} and \eqref{eq_coeff'}.
\par
The determinant equation $\text{det}(\rho)=n$ can be written as
\begin{equation}\label{eq_determinant}
c\alpha_3-f\alpha_6+j\alpha_9=n.
\end{equation}
\subsection{The core of the proof of Proposition \ref{propo_counting}}%
The proof of Proposition \ref{propo_counting} heavily relies on the following result.
\begin{proposition}\label{propo_core}
Let $x_0$, $y_0$ be some fixed integers, $D_0>0$ be an absolute constant, $U$ a large parameter, which goes to infinity and $0\leq\delta\ll 1$. Let $1\leq k\leq 5$ be an integer. Let $u$ be a real number satisfying $\abs{u}\leq U^2$, $v$ be a positive integer and $m$ be a positive integer satisfying $\abs{m}\ll U$. Let $q$ be a positive definite binary quadratic form with three uniformly bounded real coefficients of discriminant $D\geq D_0$ and $\lambda$ be a linear form on $\R^2$ with two uniformly bounded real coefficients. Assume that among the five coefficients of $q$ and $\lambda$, exactly $k$ of them are not integers. In this case,
\begin{multline*}
\left\vert\left\{(x,y)\in\Z^2, \abs{x},\abs{y}\ll U, (x,y)\equiv(x_0,y_0)\bmod{v}, q(x,y)+m\lambda(x,y)=u+O(U^2\delta)\right\}\right\vert \\
\ll_{D_0,\epsilon} U^\epsilon\left(1+\frac{U^2\delta^{1/(k+1)}}{v}\right)
\end{multline*} 
for all $\epsilon>0$. Note that the implied constant depends on $D_0$ and $\varepsilon$, but is uniform in all other parameters.
\end{proposition}
\begin{proof}[\proofname{} of Proposition \ref{propo_core}]
Let us approximate simultaneously the $k$ coefficients of $q$ and $\lambda$, say $c_1,\dots,c_k$, which are not integers by rational numbers of common denominator $1\leq r\leq R$ for some parameter $R$, which will be chosen later.
\begin{equation*}
\forall i\in\{1,\dots,k\},\quad \left\vert c_i-\frac{p_i}{r}\right\vert\leq\frac{1}{rR^{1/k}}.
\end{equation*}
If $(x,y)\in\Z^2$ satisfy $\abs{x},\abs{y}\ll U$ and
\begin{equation*}
q(x,y)+m\lambda(x,y)=u+O(U^2\delta)
\end{equation*}
then
\begin{equation*}
q_\Z(x,y)+m\lambda_\Z(x,y)=ru+O\left(RU^2\delta+\frac{U^2}{R^{1/k}}\right)
\end{equation*}
where $q_\Z$ (respectively $\lambda_\Z$) is the binary quadratic form (respectively linear form on $\R^2$) with integer coefficients obtained from $q$ (respectively $\lambda$) after substituing the coefficients $\alpha_i$ , $1\leq i\leq k$, by their rational approximation and multiplying by the common denominator $r$. The optimal choice for $R$ is given by
\begin{equation*}
R=\min{\left(\frac{U^\epsilon}{\delta^{k/(k+1)}},U^{2k}\right)}=\begin{cases}
\frac{U^\epsilon}{\delta^{k/(k+1)}} & \text{if $\delta\geq\frac{1}{U^{(k+1)(2k-\epsilon)/k}}$,} \\
U^{2k} & \text{otherwise.}
\end{cases}
\end{equation*}
In both cases, $R\to+\infty$ as $U\to\infty$ since $\delta\ll 1$. Thus, the quadratic form $r^{-1}q_\Z$, being close to the quadratic form $q$ of discriminant $D\geq D_0>0$, remains positive definite and the same holds for $q_\Z$. Note that
\begin{equation*}
q_\Z(x,y)+m\lambda_\Z(x,y)
\end{equation*}
belongs to a fixed congruence class modulo $v$. By \cite[Lemma 8 (a)]{BlPo}, the number of pairs of integers $(x,y)$ is bounded by
\begin{equation*}
\ll_\epsilon\left(RU^2+\frac{U^2}{R^{1/k}}\right)^\epsilon\left(1+\frac{RU^2\delta+\frac{U^2}{R^{1/k}}}{v}\right)\ll_\epsilon U^\epsilon\left(1+\frac{U^2\delta^{1/(k+1)}}{v}\right).
\end{equation*}
\end{proof}
\subsection{Proof of Proposition \ref{propo_counting}}%
One of the coefficients of the matrix $\rho$ is different from $0$. For instance, let us assume that $g\neq 0$ and let us set $\lambda=(g,\ell)$. There are $n^{1/3}/\lambda$ integers $g$ by \eqref{eq_coeff}. Let us fix $g$.
\par
Firstly, let us count the number of pairs $(a,d)$. The equation in \eqref{eq_column_1} can be written as
\begin{equation}\label{eq_a_d}
q_1^\R(a,d)+2g\lambda_1^\R(a,d)=n^{2/3}-(\beta_3^3-\beta_1^2x_2^2-x^2)g^2+O\left(n^{2/3}\left(\sqrt{\Delta}+\Delta\right)\right)
\end{equation}
where $q_1^\R$ is the positive definite quadratic form of discriminant $\beta_1^2\geq 3/4$ with bounded real coefficients given by
\begin{equation*}
q_1^\R(a,d)=a^2+(x_1^2+\beta_1^2)d^2-2x_1ad
\end{equation*}
and $\lambda_1^\R$ is the linear form with bounded real coefficients given by
\begin{equation*}
\lambda_1^\R(a,d)=xa-(xx_1+\beta_1^2x_2)d.
\end{equation*}
By Proposition \ref{propo_core}, the number of pairs $(a,d)$ is bounded by
\begin{equation*}
\ll_\epsilon n^\epsilon\left(1+n^{2/3}\left(\sqrt{\Delta}+\Delta\right)^{1/5}\right).
\end{equation*}
\par
Let us count the number of pairs $(h,j)$. Similarly, the equation in \eqref{eq_column_2} implies that the number of pairs $(h,j)$ is also bounded by
\begin{equation*}
\ll_\epsilon n^\epsilon\left(1+n^{2/3}\left(\sqrt{\Delta}+\Delta\right)^{1/5}\right).
\end{equation*}
\par
Let us fix $(a,d,g,h,j)$ and let us count the number of $4$-tuples $(b,c,e,f)$. We decompose this count into
\begin{equation*}
\sum_{(b,c,e,f)}1=\sum_{\substack{e \\
\ell\mid{\alpha_3}}}\sum_{\substack{f \\
\ell\mid{\alpha_2}}}\sum_{\substack{b \\
\ell\mid{\alpha_6}}}\sum_{\substack{c \\
\ell\mid{\alpha_5}}}1=\sum_{\substack{e \\
\ell\mid{\alpha_3\neq 0}}}\sum_{\substack{f \\
\ell\mid{\alpha_2}}}\sum_{\substack{b \\
\ell\mid{\alpha_6}}}\sum_{\substack{c \\
\ell\mid{\alpha_5}}}1+\sum_{\substack{e \\
\ell\mid\alpha_3=0}}\sum_{\substack{f \\
\ell\mid\alpha_2}}\sum_{\substack{c \\
\ell\mid{\alpha_5}}}\sum_{\substack{b \\
\ell\mid{\alpha_6}}}1.
\end{equation*}
Note that $\alpha_3=0$ fixes $e$. Thus, the largest count will be
\begin{equation*}
\sum_{\substack{e \\
\ell\mid{\alpha_3\neq 0}}}\sum_{\substack{f \\
\ell\mid{\alpha_2}}}\sum_{\substack{b \\
\ell\mid{\alpha_6}}}\sum_{\substack{c \\
\ell\mid{\alpha_5}}}1.
\end{equation*}
\par
Let us count the number of pairs $(e,f)$. The equation in \eqref{eq_column_3} can be writen as (after multiplying by $\beta_2^2$)
\begin{equation}\label{eq_h_j}
q_3^\R(e,f)+2\lambda_3^\R(e,f)=n^{2/3}-C+O\left(n^{2/3}\left(\sqrt{\Delta}+\Delta\right)\right)
\end{equation}
where $q_3^\R$ is the positive definite quadratic form of discriminant $\beta_2^2\geq 3/4$ with bounded real coefficients given by
\begin{equation*}
q_3^\R(e,f)=(\beta_2^2+x_2^2)e^2+f^2+2x_2ef,
\end{equation*}
$\lambda_3^\R$ is the linear form with bounded real coefficients given by
\begin{multline*}
\lambda_3^\R(e,f)=\left((\beta_2^2x_1+x_2x_3)d-(x_2^2x_3+\beta_2^2x_1x_2)g-(\beta_2^2x_2+x_2^3)h-x_2^2j\right)e \\
+\left(x_3d-x_2x_3g-x_2^2h-x_2j\right)f
\end{multline*}
and $C$ is a constant, which only depends on $z$, $d$, $g$, $h$ and $j$ and bounded by $n^{2/3}$. We will use once again Proposition \ref{propo_core} but with the additional feature that both $e$ and $f$ belong to a fixed congruence class modulo $\ell/\lambda$ since $\ell$ divides both $\alpha_2=dj-fg$ and $\alpha_3=dh-eg$. The number of pairs $(e,f)$ is bounded by
\begin{equation*}
\ll_\epsilon n^\epsilon\left(1+\frac{n^{2/3}\left(\sqrt{\Delta}+\Delta\right)^{1/5}}{\ell/\lambda}\right).
\end{equation*}
\par
Let us count the number of $b$. Equation \eqref{eq_comatrix} implies (after multiplying by $\beta_1/n^{1/3}$) that
\begin{equation*}
b=\left(\frac{\beta_1^2x_2}{n^{1/3}}g+x_1\right)e+\frac{\beta_1^2}{n^{1/3}}gf+\frac{\beta_1}{n^{1/3}}c_2+O\left(n^{1/3}\left(\sqrt{\Delta}+\Delta\right)\right)
\end{equation*}
for some constant $c_2$, which only depends on $(a,d,g,h,j)$. Moreover, $b$ belongs to a fixed congruence class modulo $\ell/\lambda$ since $\ell$ divides $\alpha_6=ah-gb$. Thus, the number of $b$ is bounded by
\begin{equation*}
1+\frac{n^{1/3}\left(\sqrt{\Delta}+\Delta\right)}{\ell/\lambda}.
\end{equation*}
\par
Let us count the number of $c$. There is only one $c$ since $c$ is fixed by the determinant equation \eqref{eq_determinant} where $\alpha_3\neq 0$. Note that this is where the condition $\alpha_3\neq 0$ is used.
\section{End of the proof of Theorem \ref{theo_mainresult}}%
\subsection{Bounding $K_{\ell,n}(z)$}%
The following proposition gives a bound for the quantities $K_{\ell,n}(z)$ given in \eqref{eq_Kmn} for any $z$ in a compact set of $X$, any positive integer $n$ and any positive integer $\ell$ dividing $n$. Let us define
\begin{equation*}
M_{\ell,n}\coloneqq\sum_{\lambda\mid\ell}\frac{1}{\lambda}\left(1+\frac{n^{2/3}}{\ell/\lambda}\right)\left(1+\frac{n^{1/3}}{\ell/\lambda}\right)
\end{equation*}
for any positive integer $n$ and any positive integer $\ell$ dividing $n$.
\begin{proposition}\label{propo_kln}
Let $n$ a positive integer, which goes to infinity with $T$ and $\ell$ a positive integer dividing $n$. If $z$ belongs to a compact subset of $X$ and $n\leq T^{3/10}$ then
\begin{equation*}
  K_{\ell,n}(z)\ll_\epsilon T^{3+\epsilon}n^{1/3+\epsilon}+T^{2+\epsilon}n^{5+\epsilon}M_{\ell,n}.
\end{equation*}
\end{proposition}
\begin{proof}[\proofname{} of Proposition \ref{propo_kln}]
By Proposition \ref{propo_estimates}, if $1\leq\alpha_1(a), \alpha_2(a)\ll 1$ then
\begin{equation*}
\mathcal{H}^{-1}\left(h_T^W\right)(a)\ll\begin{cases}
  T^{3+\epsilon} & \text{if $1\leq\alpha_1(a),\alpha_2(a)\leq1+1/n^{10/3}$,} \\
  T^{3/2+\epsilon}n^5 & \text{$1+1/n^{10/3}\leq\alpha_1(a),\alpha_2(a)\ll 1$,} \\
  T^{2+\epsilon} n^{10/3} & \text{otherwise.}
\end{cases}
\end{equation*}
By Proposition \ref{propo_counting}, if $0\leq\delta_1, \delta_2\ll 1$ then
\begin{equation*}
M_{\ell,n}(z;\delta_1,\delta_2)\ll\begin{cases}
n^{1/3+\epsilon} & \text{if $0\leq\delta_1,\delta_2\leq 1/n^{10/3}$,} \\
n^{5/3+\epsilon}M_{\ell,n} & \text{otherwise.}
\end{cases}
\end{equation*}
These two facts conclude the proof since if $n\leq T^{3/10}$ then $T^{3/2}n^{20/3}\leq T^2n^5$.
\end{proof}
\subsection{Proof of Theorem \ref{theo_mainresult}}%
Let us quickly finish the proof of Theorem \ref{theo_mainresult}. By Rankin-Selberg theory and the Cauchy-Schwarz inequality, the amplifier defined in \eqref{eq_choice_ampli} satisfies
\begin{eqnarray*}
\vert\vert\alpha\vert\vert_2^2 & \ll_\epsilon & L^{1+\epsilon}, \\
\vert\vert\alpha\vert\vert_1 & \ll_\epsilon & L^{1+\epsilon}
\end{eqnarray*}
for any $\epsilon>0$.
\par
Thus, by Proposition \ref{propo_first_11} and Proposition \ref{propo_kln}, if $L\leq T^{3/40}$ then
\begin{equation*}
  \left\vert\Phi_{j_0}(z)\right\vert^2\ll_\epsilon (TL)^\epsilon\left(\frac{T^3}{L}+T^2L^{18}\right).
\end{equation*}
The optimal choice for $L$ is given by $L=T^{1/19}\leq T^{3/40}$, which implies Theorem \ref{theo_mainresult}.
\bibliographystyle{plain}
\bibliography{biblioo}
\end{document}